\documentclass[12pt]{amsart}      
\usepackage{amssymb}
\usepackage{eucal}
\usepackage{amsmath}
\usepackage{amscd}
\usepackage[all]{xy}           

\newcommand{\nc}{\newcommand}
\newcommand{\delete}[1]{}
\nc{\dfootnote}[1]{{}}          
\nc{\ffootnote}[1]{\dfootnote{#1}}
\nc{\mfootnote}[1]{\footnote{#1}} 
\nc{\ofootnote}[1]{\footnote{\tiny Older version: #1}} 
\nc{\mlabel}[1]{\label{#1}}  
\nc{\mcite}[1]{\cite{#1}  
{\hfill \hspace{1cm}{\bf{{\ }\hfill(#1)}}}}
\nc{\mkeep}[1]{\marginpar{{\bf #1}}} 
\nc{\newpart}[1]{{\bf NEW: } {\bf #1}} 

\newtheorem{theorem}{Theorem}[section]
\newtheorem{prop}[theorem]{Proposition}
\newtheorem{defn}[theorem]{Definition}
\newtheorem{lemma}[theorem]{Lemma}
\newtheorem{coro}[theorem]{Corollary}
\newtheorem{prop-def}{Proposition-Definition}[section]

\nc{\dfgen}{\Omega}
\nc{\dfrel}{\Lambda}
\nc{\dfgenb}{\vec{\omega}}
\nc{\dfrelb}{\vec{\lambda}}
\nc{\dfgene}{\omega}
\nc{\dfrele}{\lambda}
\nc{\dftimes}{\,\square\,}      
\nc{\dfotimes}{\,\maltese\,}    
\DeclareMathOperator*{\dfTimes}{\Large{\square}}
\nc{\dfl}{\succ}
\nc{\dfr}{\prec}
\nc{\dfc}{\circ}
\nc{\dfb}{\bullet}
\nc{\dft}{\star}
\nc{\dfpair}[2]{\left[ {{#1}\atop{#2}}\right]}
\nc{\dfll}{{\dfpair{\dfl}{\dfl}}}
\nc{\dflr}{{\dfpair{\dfl}{\dfr}}}
\nc{\dflc}{{\dfpair{\dfl}{\dfc}}}
\nc{\dfrl}{{\dfpair{\dfr}{\dfl}}}
\nc{\dfrr}{{\dfpair{\dfr}{\dfr}}}
\nc{\dfrc}{{\dfpair{\dfr}{\dfc}}}
\nc{\dfcl}{{\dfpair{\dfc}{\dfl}}}
\nc{\dfcr}{{\dfpair{\dfc}{\dfr}}}
\nc{\dfcc}{{\dfpair{\dfc}{\dfc}}}
\nc{\dflt}{{\dfpair{\dfl}{\dft}}}
\nc{\dfrt}{{\dfpair{\dfr}{\dft}}}
\nc{\dfct}{{\dfpair{\dfc}{\dft}}}
\nc{\dftl}{{\dfpair{\dft}{\dfl}}}
\nc{\dftr}{{\dfpair{\dft}{\dfr}}}
\nc{\dftc}{{\dfpair{\dft}{\dfc}}}
\nc{\dftt}{{\dfpair{\dft}{\dft}}}
\nc{\dflb}{{\dfpair{\dfl}{\bullet}}}
\nc{\dfrb}{{\dfpair{\dfr}{\bullet}}}
\nc{\dftb}{{\dfpair{\dft}{\bullet}}}
\nc{\dfbt}{{\dfpair{\bullet}{\dft}}}
\nc{\dfbb}{{\dfpair{\bullet}{\bullet}}}
\nc{\dfbl}{{\dfpair{\bullet}{\dfl}}}
\nc{\dfbr}{{\dfpair{\bullet}{\dfr}}}
\nc{\dfop}{\odot}
\nc{\dfoa}{\dfop^{(1)}}
\nc{\dfob}{\dfop^{(2)}}
\nc{\dfoc}{\dfop^{(3)}}
\nc{\dfod}{\dfop^{(4)}}
\nc{\tstar}{\tilde{\star}}
\nc{\tdfop}{\tilde{\dfop}}
\nc{\tdfoa}{\tilde{\dfop}^{(1)}}
\nc{\tdfob}{\tilde{\dfop}^{(2)}}
\nc{\tdfoc}{\tilde{\dfop}^{(3)}}
\nc{\tdfod}{\tilde{\dfop}^{(4)}}
\nc{\dflls}{{{\dfl}\atop{\dfl}}}
\nc{\dflrs}{{{\dfl}\atop{\dfr}}}
\nc{\dflcs}{{{\dfl}\atop{\dfc}}}
\nc{\dfrls}{{{\dfr}\atop{\dfl}}}
\nc{\dfrrs}{{{\dfr}\atop{\dfr}}}
\nc{\dfrcs}{{{\dfr}\atop{\dfc}}}
\nc{\dfcls}{{{\dfc}\atop{\dfl}}}
\nc{\dfcrs}{{{\dfc}\atop{\dfr}}}
\nc{\dfccs}{{{\dfc}\atop{\dfc}}}
\nc{\dflts}{{{\dfl}\atop{\dft}}}
\nc{\dfrts}{{{\dfr}\atop{\dft}}}
\nc{\dfcts}{{{\dfc}\atop{\dft}}}
\nc{\dftls}{{{\dft}\atop{\dfl}}}
\nc{\dftrs}{{{\dft}\atop{\dfr}}}
\nc{\dftcs}{{{\dft}\atop{\dfc}}}
\nc{\dftts}{{{\dft}\atop{\dft}}}
\nc{\ndfll}{{\dfpair{\dfl}{>}}}
\nc{\ndflr}{{\dfpair{\dfl}{<}}}
\nc{\ndflc}{{\dfpair{\dfl}{\bullet}}}
\nc{\ndfrl}{{\dfpair{\dfr}{>}}}
\nc{\ndfrr}{{\dfpair{\dfr}{<}}}
\nc{\ndfrc}{{\dfpair{\dfr}{\bullet}}}
\nc{\ndfcl}{{\dfpair{\dfc}{>}}}
\nc{\ndfcr}{{\dfpair{\dfc}{<}}}
\nc{\ndfcc}{{\dfpair{\dfc}{\bullet}}}
\nc{\ndflt}{{\dfpair{\dfl}{\diamondsuit}}}
\nc{\ndfrt}{{\dfpair{\dfr}{\diamondsuit}}}
\nc{\ndfct}{{\dfpair{\dfc}{\diamondsuit}}}
\nc{\ndftl}{{\dfpair{\dft}{>}}}
\nc{\ndftr}{{\dfpair{\dft}{<}}}
\nc{\ndftc}{{\dfpair{\dft}{\bullet}}}
\nc{\ndftt}{{\dfpair{\dft}{\diamondsuit}}}
\nc{\dftri}[3]{{\left[{{#1}\atop{#2}}\atop{\scriptstyle{#3}}\right]}}
\nc{\dfrrr}{{\dftri{\dfr}{\dfr}{\dfr}}}
\nc{\dfrrl}{{\dftri{\dfr}{\dfr}{\dfl}}}
\nc{\dfrrb}{{\dftri{\dfr}{\dfr}{\dfb}}}
\nc{\dfrrt}{{\dftri{\dfr}{\dfr}{\dft}}}
\nc{\dfrlr}{{\dftri{\dfr}{\dfl}{\dfr}}}
\nc{\dfrll}{{\dftri{\dfr}{\dfl}{\dfl}}}
\nc{\dfrlb}{{\dftri{\dfr}{\dfl}{\dfb}}}
\nc{\dfrlt}{{\dftri{\dfr}{\dfl}{\dft}}}

\nc{\dfrbr}{{\dftri{\dfr}{\dfb}{\dfr}}}
\nc{\dfrbl}{{\dftri{\dfr}{\dfb}{\dfl}}}
\nc{\dfrbb}{{\dftri{\dfr}{\dfb}{\dfb}}}
\nc{\dfrbt}{{\dftri{\dfr}{\dfb}{\dft}}}
\nc{\dfrtr}{{\dftri{\dfr}{\dft}{\dfr}}}
\nc{\dfrtl}{{\dftri{\dfr}{\dft}{\dfl}}}
\nc{\dfrtb}{{\dftri{\dfr}{\dft}{\dfb}}}
\nc{\dfrtt}{{\dftri{\dfr}{\dft}{\dft}}}

\nc{\dflrr}{{\dftri{\dfl}{\dfr}{\dfr}}}
\nc{\dflrl}{{\dftri{\dfl}{\dfr}{\dfl}}}
\nc{\dflrb}{{\dftri{\dfl}{\dfr}{\dfb}}}
\nc{\dflrt}{{\dftri{\dfl}{\dfr}{\dft}}}
\nc{\dfllr}{{\dftri{\dfl}{\dfl}{\dfr}}}
\nc{\dflll}{{\dftri{\dfl}{\dfl}{\dfl}}}
\nc{\dfllb}{{\dftri{\dfl}{\dfl}{\dfb}}}
\nc{\dfllt}{{\dftri{\dfl}{\dfl}{\dft}}}

\nc{\dflbr}{{\dftri{\dfl}{\dfr}{\dfr}}}
\nc{\dflbl}{{\dftri{\dfl}{\dfr}{\dfl}}}
\nc{\dflbb}{{\dftri{\dfl}{\dfr}{\dfb}}}
\nc{\dflbt}{{\dftri{\dfl}{\dfr}{\dft}}}
\nc{\dfltr}{{\dftri{\dfl}{\dft}{\dfr}}}
\nc{\dfltl}{{\dftri{\dfl}{\dft}{\dfl}}}
\nc{\dfltb}{{\dftri{\dfl}{\dfl}{\dfb}}}
\nc{\dfltt}{{\dftri{\dfl}{\dft}{\dft}}}

\nc{\dfbrr}{{\dftri{\dfb}{\dfr}{\dfr}}}
\nc{\dfbrl}{{\dftri{\dfb}{\dfr}{\dfl}}}
\nc{\dfbrb}{{\dftri{\dfb}{\dfr}{\dfb}}}
\nc{\dfbrt}{{\dftri{\dfb}{\dfr}{\dft}}}
\nc{\dfblr}{{\dftri{\dfb}{\dfl}{\dfr}}}
\nc{\dfbll}{{\dftri{\dfb}{\dfl}{\dfl}}}
\nc{\dfblb}{{\dftri{\dfb}{\dfl}{\dfb}}}
\nc{\dfblt}{{\dftri{\dfb}{\dfl}{\dft}}}

\nc{\dfbbr}{{\dftri{\dfb}{\dfb}{\dfr}}}
\nc{\dfbbl}{{\dftri{\dfb}{\dfb}{\dfl}}}
\nc{\dfbbb}{{\dftri{\dfb}{\dfb}{\dfb}}}
\nc{\dfbbt}{{\dftri{\dfb}{\dfb}{\dft}}}
\nc{\dfbtr}{{\dftri{\dfb}{\dft}{\dfr}}}
\nc{\dfbtl}{{\dftri{\dfb}{\dft}{\dfl}}}
\nc{\dfbtb}{{\dftri{\dfb}{\dft}{\dfb}}}
\nc{\dfbtt}{{\dftri{\dfb}{\dft}{\dft}}}

\nc{\dftrr}{{\dftri{\dft}{\dfr}{\dfr}}}
\nc{\dftrl}{{\dftri{\dft}{\dfr}{\dfl}}}
\nc{\dftrb}{{\dftri{\dft}{\dfr}{\dfb}}}
\nc{\dftrt}{{\dftri{\dft}{\dfr}{\dft}}}
\nc{\dftlr}{{\dftri{\dft}{\dfl}{\dfr}}}
\nc{\dftll}{{\dftri{\dft}{\dfl}{\dfl}}}
\nc{\dftlb}{{\dftri{\dft}{\dfl}{\dfb}}}
\nc{\dftlt}{{\dftri{\dft}{\dfl}{\dft}}}

\nc{\dftbr}{{\dftri{\dft}{\dfb}{\dfr}}}
\nc{\dftbl}{{\dftri{\dft}{\dfb}{\dfl}}}
\nc{\dftbb}{{\dftri{\dft}{\dfb}{\dfb}}}
\nc{\dftbt}{{\dftri{\dft}{\dfb}{\dft}}}
\nc{\dfttr}{{\dftri{\dft}{\dft}{\dfr}}}
\nc{\dfttl}{{\dftri{\dft}{\dft}{\dfl}}}
\nc{\dfttb}{{\dftri{\dft}{\dft}{\dfb}}}
\nc{\dfttt}{{\dftri{\dft}{\dft}{\dft}}}

\renewcommand{\ne}{\nearrow} 
\newcommand{\se}{\searrow}
\newcommand{\sw}{\swarrow}
\newcommand{\nw}{\nwarrow}

\newcommand{\east}{\succ}
\newcommand{\west}{\prec}
\newcommand{\north}{\wedge}
\newcommand{\south}{\vee}

\nc{\bin}[2]{ (_{\stackrel{\scs{#1}}{\scs{#2}}})}  
\nc{\binc}[2]{ \left (\!\! \begin{array}{c} \scs{#1}\\
    \scs{#2} \end{array}\!\! \right )}  
\nc{\bincc}[2]{  \left ( {\scs{#1} \atop
    \vspace{-.5cm}\scs{#2}} \right )}  
\nc{\bs}{\bar{S}}
\nc{\la}{\longrightarrow}
\nc{\rar}{\rightarrow}
\nc{\dar}{\downarrow}
\nc{\dap}[1]{\downarrow \rlap{$\scriptstyle{#1}$}}
\nc{\uap}[1]{\uparrow \rlap{$\scriptstyle{#1}$}}
\nc{\defeq}{\stackrel{\rm def}{=}}
\nc{\oeq}[1]{\stackrel{(#1)}{=}}
\nc{\dis}[1]{\displaystyle{#1}}
\nc{\dotcup}{\ \displaystyle{\bigcup^\bullet}\ }
\nc{\hcm}{\ \hat{,}\ }
\nc{\hcirc}{\hat{\circ}}
\nc{\hts}{\hat{\otimes}}
\nc{\lts}{\stackrel{\leftarrow}{\otimes}}
\nc{\rts}{\stackrel{\rightarrow}{\otimes}}
\nc{\lleft}{[}
\nc{\lright}{]}
\nc{\curlyl}{\left \{ \begin{array}{c} {} \\ {} \end{array}
    \right .  \!\!\!\!\!\!\!}
\nc{\curlyr}{ \!\!\!\!\!\!\!
    \left . \begin{array}{c} {} \\ {} \end{array}
    \right \} }
\nc{\longmid}{\left | \begin{array}{c} {} \\ {} \end{array}
    \right . \!\!\!\!\!\!\!}
\nc{\ora}[1]{\stackrel{#1}{\rar}}
\nc{\ola}[1]{\stackrel{#1}{\la}}
\nc{\olrarrow}[1]{\stackrel{#1}{\Leftrightarrow}}
\nc{\scs}[1]{\scriptstyle{#1}}
\nc{\mrm}[1]{{\rm #1}}
\nc{\margin}[1]{\marginpar{\rm #1}}   
\nc{\dirlim}{\displaystyle{\lim_{\longrightarrow}}\,}
\nc{\invlim}{\displaystyle{\lim_{\longleftarrow}}\,}
\nc{\mvp}{\vspace{0.5cm}}
\nc{\ot}{\otimes}
\nc{\tk}{^{(k)}}
\nc{\tp}{^\prime}
\nc{\ttp}{^{\prime\prime}}
\nc{\svp}{\vspace{2cm}}
\nc{\vp}{\vspace{8cm}}
\nc{\proofbegin}{\noindent{\bf Proof: }}
\nc{\proofend}{$\blacksquare$ \vspace{0.5cm}}
\nc{\modg}[1]{\!<\!\!{#1}\!\!>}
\nc{\intg}[1]{F_C(#1)}
\nc{\lmodg}{\!<\!\!}
\nc{\rmodg}{\!\!>\!}
\nc{\cpi}{\widehat{\Pi}}
\nc{\sha}{{\mbox{\cyr X}}}  
\nc{\shpr}{\diamond}    
\nc{\vep}{\varepsilon}
\nc{\labs}{\mid\!}
\nc{\rabs}{\!\mid}
\nc{\hsha}{\widehat{\sha}}
\nc{\lsha}{\stackrel{\leftarrow}{\sha}}
\nc{\rsha}{\stackrel{\rightarrow}{\sha}}

\nc{\ann}{\mrm{ann}}
\nc{\Aut}{\mrm{Aut}}
\nc{\can}{\mrm{can}}
\nc{\colim}{\mrm{colim}}
\nc{\Cont}{\mrm{Cont}}
\nc{\rchar}{\mrm{char}}
\nc{\cok}{\mrm{coker}}
\nc{\dtf}{{R-{\rm tf}}}
\nc{\dtor}{{R-{\rm tor}}}

\nc{\Div}{{\mrm Div}}
\nc{\End}{\mrm{End}}
\nc{\Ext}{\mrm{Ext}}
\nc{\Fil}{\mrm{Fil}}
\nc{\Frob}{\mrm{Frob}}
\nc{\Gal}{\mrm{Gal}}
\nc{\GL}{\mrm{GL}}
\nc{\Hom}{\mrm{Hom}}
\nc{\hsr}{\mrm{H}}
\nc{\hpol}{\mrm{HP}}
\nc{\id}{\mrm{id}}
\nc{\im}{\mrm{im}}
\nc{\incl}{\mrm{incl}}
\nc{\length}{\mrm{length}}
\nc{\Loday}{ABQR }
\nc{\mchar}{\rm char}
\nc{\mpart}{\mrm{part}}
\nc{\ql}{{\QQ_\ell}}
\nc{\qp}{{\QQ_p}}
\nc{\rank}{\mrm{rank}}
\nc{\rcot}{\mrm{cot}}
\nc{\rdef}{\mrm{def}}
\nc{\rdiv}{{\rm div}}
\nc{\rtf}{{\rm tf}}
\nc{\rtor}{{\rm tor}}
\nc{\res}{\mrm{res}}
\nc{\SL}{\mrm{SL}}
\nc{\Spec}{\mrm{Spec}}
\nc{\tor}{\mrm{tor}}
\nc{\Tr}{\mrm{Tr}}
\nc{\tr}{\mrm{tr}}

\nc{\ab}{\mathbf{Ab}}
\nc{\Alg}{\mathbf{Alg}}
\nc{\Bax}{\mathbf{Bax}}
\nc{\bfk}{{\bf k}}
\nc{\bfone}{{\bf 1}}
\nc{\base}[1]{{a_{#1}}}
\nc{\detail}{\marginpar{\bf More detail}
    \noindent{\bf Need more detail!}
    \svp}
\nc{\Diff}{\mathbf{Diff}}
\nc{\gap}{\marginpar{\bf Incomplete}\noindent{\bf Incomplete!!}
    \svp}
\nc{\FMod}{\mathbf{FMod}}
\nc{\Int}{\mathbf{Int}}
\nc{\Mon}{\mathbf{Mon}}
\nc{\remarks}{\noindent{\bf Remarks: }}
\nc{\Rep}{\mathbf{Rep}}
\nc{\Rings}{\mathbf{Rings}}
\nc{\Sets}{\mathbf{Sets}}
\nc{\bill}[1]{\marginpar{\bf To Bill}\noindent{\bf To Bill:}
    {\tt #1}\\ }
\nc{\li}[1]{\marginpar{\bf To Li}\noindent{\bf To Li:}
    {\tt #1}\\ }

\nc{\BA}{{\Bbb A}}
\nc{\CC}{{\Bbb C}}
\nc{\DD}{{\Bbb D}}
\nc{\EE}{{\Bbb E}}
\nc{\FF}{{\Bbb F}}
\nc{\GG}{{\Bbb G}}
\nc{\HH}{{\Bbb H}}
\nc{\LL}{{\Bbb L}}
\nc{\NN}{{\Bbb N}}
\nc{\QQ}{{\Bbb Q}}
\nc{\RR}{{\Bbb R}}
\nc{\TT}{{\Bbb T}}
\nc{\VV}{{\Bbb V}}
\nc{\ZZ}{{\Bbb Z}}

\nc{\cala}{{\mathcal A}}
\nc{\calc}{{\mathcal C}}
\nc{\cald}{{\mathcal D}}
\nc{\cale}{{\mathcal E}}
\nc{\calf}{{\mathcal F}}
\nc{\calg}{{\mathcal G}}
\nc{\calh}{{\mathcal H}}
\nc{\cali}{{\mathcal I}}
\nc{\call}{{\mathcal L}}
\nc{\calm}{{\mathcal M}}
\nc{\caln}{{\mathcal N}}
\nc{\calo}{{\mathcal O}}
\nc{\calp}{{\mathcal P}}
\nc{\calq}{{\mathcal Q}}
\nc{\calr}{{\mathcal R}}
\nc{\calt}{{\mathcal T}}
\nc{\calw}{{\mathcal W}}
\nc{\calx}{{\mathcal X}}
\nc{\CA}{\mathcal{A}}

\nc{\fraka}{{\frak a}}
\nc{\frakB}{{\frak B}}
\nc{\frakm}{{\frak m}}
\nc{\frakp}{{\frak p}}

\font\cyr=wncyr10

\begin{document}

\title{On Products and Duality of Binary, Quadratic, Regular Operads
}
\author{KURUSCH EBRAHIMI-FARD}
\address{Universit\"at Bonn -
         Physikalisches Institut,
         Nussallee 12,
         D-53115 Bonn,
         Germany}
\email{kurusch@ihes.fr}
\author{Li Guo}
\address{
Department of Mathematics and Computer Science,
Rutgers University,
Newark, NJ 07102, USA}
\email{liguo@newark.rutgers.edu}

\date{\today}


\begin{abstract}
Since its introduction by Loday in 1995 with motivation from
algebraic $K$-theory, dendriform dialgebras have been studied
quite extensively with connections to several areas in mathematics
and physics. A few more similar structures have been found
recently, such as the tri-, quadri-, ennea- and octo-algebras,
with increasing complexity in their constructions and properties.
We consider these constructions as operads and their products and duals, 
in terms of generators and relations, with the goal to clarify and
simplify the process of obtaining new algebra structures from
known structures and from linear operators. 
\end{abstract}

\maketitle




\setcounter{section}{0}

\section{Introduction}
\label{intro}

In order to study the periodicity of algebraic $K$-groups,
J.-L. Loday laid out a program in~\cite{Lo1} which led him to
the concepts of associative dialgebra and dendriform
dialgebra~\cite{Lo2}. In the next few years, their properties were
studied by several authors in areas related to operads~\cite{Lo6},
homology~\cite{Fra1,Fra2}, Hopf algebras~\cite{Ch,Hol,Ron,L-R3},
combinatorics~\cite{Fo,L-R1,A-S1,A-S2}, arithmetic~\cite{Lo5} and
quantum field theory~\cite{Fo}. See~\cite{Lo4} and other articles
in the volume for a survey of some of these
developments.

Since 2002, quite a few more similar algebra structures have been
introduced, such as the associative trialgebra and dendriform
trialgebra of Loday and Ronco~\cite{L-R2}, the dendriform
quadri-algebra of Aguiar and Loday~\cite{A-L}, the ennea-algebra,
the NS-algebra, the dendriform-Nijenhuis algebra and the
octo-algebra of Leroux~\cite{Le1,Le2,Le3}. These algebras have a
common property of ``splitting associativity", that is, expressing
the multiplication of an associative algebra as the sum of a
string of binary operations. The operations in the string satisfy
a set of relations and the associativity of the multiplication
follows from the sum of these relations. The first instance of
such algebras, the dendriform dialgebra, has a string of two
operators. The later constructions were largely inspired by the
connection~\cite{Ag,EF1,Le1} with Rota-Baxter operators\footnote{They used to be
called Baxter operators. They are renamed Rota-Baxter operators to
distinguish it clearly from the very related Yang-Baxter
operators. The latter Baxter is the Australian physicist Rodney
Baxter.} which were
introduced by G. Baxter~\cite{Ba} in 1960 and were actively
studied in the 1970s~\cite{Ro1,Ro2} and again in recent years in
connection with several areas of mathematics and
physics~\cite{A-G-K-O,C-K1,C-K2,EF2,E-G,E-G-K2,E-G-K2,Gu-Ke1,Gu-Ke2,
Gu2,Gu4,Ho}.

Two themes can be found in these recent constructions. One is the
construction of a new type of algebras that has the combined
features of types of two or more algebras that were previously
known. The other is the use of a linear operator with certain
features, such as a Rota-Baxter operator, on a known type of algebras 
to obtain another type of algebras with richer structures. 
Even though the ideas of the themes are simple, to carry them out for 
a particular construction can be quite complicated. 

The purpose of this paper is to study these constructions in the
framework of operads and their products, given by generators and
relations. This enables us to clarify, simplify and further
generalize the constructions and properties of these recent
algebra structures. 

Here is a more detailed plan of the paper. In ~\S\ref{sec:def}, we
recall the operads that give rise to the above algebra
structures. These operads are binary, quadratic and regular~\cite{Lo6}
operads with a splitting of associativity. 
To ease the notation, we call them \Loday operads and the corresponding algebras 
\Loday algebras. 
The generator-relation
construction of \Loday operads allows a concrete description which
is quite simple and can be found in the existing
literature~\cite{Lo4}. We use this description to formally define
the types of such operads.

We then define in \S\ref{sec:prod} products of \Loday algebras that 
are similar to (but different from) the operad products of 
Manin-Ginzburg-Kapranov~\cite{Ma1,Ma2,G-K,Lo3}. We show, formalizing the 
first theme, that some recently obtained \Loday algebras~\cite{A-L,Le1,Le2,Le3} 
are products of simpler algebras. Properties of
products of \Loday algebras are also studied in this section. 
The subsequent Section \ref{sec:dual} considers the dual of a \Loday operad 
and its relation with the products.  

A full understanding of the second theme mentioned above depends on sufficient
knowledge of the linear operators that give rise to \Loday algebras. While this
topic is being investigated in another project, using the framework introduced here,
we can make the second theme precise for most of the operators that we are aware of,
when the operators are applied to any types of \Loday algebras. This is presented
in \S\ref{sec:op}.

The concept of unit actions of operads has recently been introduced by 
Loday~\cite{Lo6} and used to construct Hopf algebras on the free algebras. 
In a separate work~\cite{E-G2}, we investigate the relation between products 
of operads and their unit actions. See~\cite{Va} also for a more recent 
application of products of \Loday operads. 

\section{\Loday algebras and operads}
\mlabel{sec:def}

\subsection{\Loday algebras and their types}
Let $\bfk$ be a field of characteristic zero. Let $A$ be a vector space over $\bfk$.
The dendriform dialgebra and its generalizations, such as the trialgebra,
quadri-algebra, ennea-algebra or Nijenhuis-dendriform algebra,
are constructed from a finite set of binary operations
$$\dfgen:=\{\dfop_n: A\otimes_\bfk A\to A,\, n=1,\cdots,m\}$$
together with a set of ``associativity" relations. This
construction fits into the general framework of operads~\cite{Fre,G-K,Lo3,Lo6} 
that are binary, quadratic and regular. We will first briefly recall 
these concepts and refer to the above references for details. 
We then give a more concrete
definition of such algebras adapted from~\cite{Lo4,L-R2}.
This definition is easier to work with and applies when the base
field $\bfk$ is replaced by a commutative ring  with identity.


An (algebraic) operad is a sequence $\{\calp(n)\}$ of finitely generated
$\bfk[S_n]$-modules such that the Schur functor 
$$ \calp: V\mapsto \bigoplus_n (\calp(n) \ot V^{\ot n})_{S_n}$$
is equipped with a composition law $\gamma: \calp\circ \calp \to \calp$ which 
is associative and unital. 
An operad $\{\calp(n)\}$ is called {\bf binary} if $\calp(1)=\bfk$ and $\calp(n),
n\geq 3$ are induced from $\calp(2)$ by composition; is called {\bf quadratic} if
all relations among the binary operations in $\calp(2)$ are derived
from $\calp(3)$; is called {\bf regular} if,
moreover, the binary operations have no symmetries (such as
$x\cdot y=y\cdot x$), and the induced relations in $\calp(3)$
occur in the same order (such as $(x\cdot y) \cdot z=x\cdot
(y\cdot z)$, not $(x\cdot y)\cdot z=x\cdot (z\cdot y)$). 

By regularity, the space $\calp(n)$ is of the form $\calp_n\otimes
\bfk[S_n]$ where $\calp_n$ is a vector space. So the operad $\{\calp(n)\}$ is
determined by $\{\calp_n\}$. Then a binary, quadratic, regular operad 
is determined by a pair $(\dfgen,\dfrel)$ where $\dfgen=\calp_2$, called 
the {\bf space of generators},  
and $\dfrel$ is a subspace of $\dfgen^{\otimes 2} \oplus \dfgen^{\otimes 2}$, 
called the {\bf space of relations.} The pair $(\dfgen,\dfrel)$ is
called the {\bf type} of the operad $\calp$, or of the
corresponding algebra structure.

A $\bfk$-vector space $A$ is called a {\bf binary, quadratic, regular algebra 
of type $(\dfgen,\dfrel)$} if it has binary operations $\dfgen$ and if,
for 
$$\big (\sum_{i=1}^k\dfoa_i\otimes \dfob_i, \sum_{j=1}^m\dfoc_j\otimes \dfod_j \big)
\in \dfrel \subseteq
\dfgen^{\otimes 2} \oplus \dfgen^{\otimes 2}$$ with
$\dfoa_i,\dfob_i,\dfoc_j,\dfod_j\in \dfgen$, $1\leq i\leq k$, $1\leq j\leq m$, 
we have
\begin{equation}
 \sum_{i=1}^k(x\dfoa_i y) \dfob_i z = \sum_{j=1}^m x \dfoc_j (y \dfod_j z), \forall\ x,y,z\in A.
\mlabel{eq:rel}
\end{equation}
When there is no danger of confusion, we will use $\dfrel$
to denote the type of an algebra structure.

Since a type $(\dfgen,\dfrel)$ is determined by
$(\dfgenb,\dfrelb)$ where $\dfgenb$ is a basis of $\dfgen$ and
$\dfrelb$ is a basis of $\dfrel$, we also use $(\dfgenb,\dfrelb)$
to denote for the type of a binary, quadratic, regular operad, as is 
usually the case in the literature.

We say that a binary, quadratic, regular operad $(\dfgen,\dfrel)$ 
{\bf has a splitting associativity} if there is an element $\star$ of 
$\dfgen$ such that $(\star\otimes \star, \star\otimes \star)$ is 
in $\dfrel$~\cite{Lo6}. 
As abbreviation, we call such an operad an {\bf associative BQR operad}, or 
simply an {\bf \Loday operad}. 

\begin{lemma}
A binary, quadratic, regular operad is \Loday if and only if there is a basis
$\dfgenb=\{\dfgene_i\}$ of $\dfgen$ such that $\star=\sum_i
\dfgene_i$ and there is a
basis $\dfrelb=\{\dfrele_j\}_j$ of $\dfrel$ such that the
associativity of $\star$ is given by the sum of $\dfrele_j$ (splitting
associativity): 
$$ (\star\otimes \star, \star\otimes \star)= \sum_j \dfrele_j.$$
\end{lemma}
\begin{proof}
The if part is clear. For the only if part, given $\star\in
\dfgen$ such that $(\star\otimes \star, \star\otimes \star)$ is in
$\dfrel$, complete $\star$ to a basis
$\{\star,\dfgene_2,\cdots,\dfgene_r\}$ of $\dfgen$ and then take
$\dfgene_1=\star-\dfgene_2-\cdots - \dfgene_r$. Similarly
complete $(\star\otimes \star, \star\otimes \star)$ to a basis
$\{(\star\otimes\star,\star\otimes \star)
,\dfrele_2,\cdots,\dfrele_s\}$ of $\dfrel$ and then take
$\dfrele_1=(\star\otimes \star,\star\otimes
\star)-\dfrele_2-\cdots - \dfrele_s$.
\end{proof}

Let $(\dfgen,\dfrel)$ and $(\dfgen',\dfrel')$ be \Loday operads with 
associative operations $\star$ and $\star'$ respectively. 
A {\bf morphism} $f:(\dfgen,\dfrel)\to (\dfgen',\dfrel')$ is a linear map
$\dfgen\to \dfgen'$ sending $\star$ to $\star'$ and inducing a linear map 
$\dfrel\to \dfrel'$. An invertible morphism is called an {\bf isomorphism}, 
and called an {\bf automorphism} if $(\dfgen,\dfrel)=(\dfgen',\dfrel')$.

\subsection{Examples}
We now describe known examples of dendriform related algebras in the
context of operad types that were just defined, as illustrations
of the concepts and as preparations for later applications.

\smallskip
\noindent 1. (Associative algebra) An associative $\bfk$-algebra
is a $\bfk$-vector space $A$ with an associative \mbox{product
$\cdot$\,.} This means that it is of type $(\dfgenb_A,\dfrelb_A)$
with $\dfgenb_A=\{\cdot\}$ and
$\dfrelb_A=\{(\cdot\otimes\cdot,\cdot\otimes\cdot)\}.$
\smallskip\\
2. (Dialgebra) The {\bf dendriform dialgebra} of Loday~\cite{Lo4} is
defined to be a $\bfk$-vector space $D$ with binary operations
$\prec$ and  $\succ$ such that
$$(x\prec y)\prec z=x\prec (y\prec z+y\succ z),
(x\succ y)\prec z=x\succ (y\prec z), $$
$$(x\prec y+x\succ y)\succ z=x\succ (y\succ z).$$
This means that $D$ is of operad type $(\dfgenb_D,\dfrelb_D)$ with
$\dfgenb_D=\{\prec,\succ\}$ and
\begin{equation}
\dfrelb_D=\{(\prec\otimes\prec,\prec\otimes(\prec+\succ)),
(\succ\otimes\prec,\succ\otimes\prec),
((\prec+\succ)\otimes\succ,\succ\otimes\succ)\}.
\mlabel{eq:dia}
\end{equation}
The associativity of $\star:=\prec+\succ$ follows from the sum of these
equations.
If we exchange $\prec$ and $\succ$ in both $\dfgenb_D$ and $\dfrelb_D$,
we get a type $(\dfgenb^{op}_D, \dfrelb^{op}_D)$ that is isomorphic
to $(\dfgenb_D,\dfrelb_D)$. It is called the opposite of
$(\dfgenb_D,\dfrelb_D)$. The same holds for
the dendriform trialgebra, and Nijenhuis trialgebra in the following examples.
\smallskip \\
3. (Trialgebra) The {\bf dendriform trialgebra} of Loday and
Ronco~\cite{L-R2} is a $\bfk$-vector space $T$ equipped with
binary operations $\prec,\succ$ and $\circ$ that satisfy the
relations
\begin{eqnarray*}
&&(x\prec y)\prec z=x\prec (y\star z),
(x\succ y)\prec z=x\succ (y\prec z),\\
&&(x\star y)\succ z=x\succ (y\succ z),
(x\succ y)\circ z=x\succ (y\circ z), \\
&&(x\prec y)\circ z=x\circ (y\succ z),
(x\circ y)\prec z=x\circ (y\prec z),
(x\circ y)\circ z=x\circ (y\circ z)
\end{eqnarray*}
for $x,y,z\in D$. Here $\star=\prec+\succ+\circ.$
This is the \Loday algebra of type
$(\dfgenb_T,\dfrelb_T)$ with
$\dfgenb_T=\{\prec,\succ,\circ\}$ and
\begin{eqnarray}
\dfrelb_T=& \{(\prec\otimes\prec,\prec\otimes\star),
(\succ\otimes\prec,\succ\otimes\prec),
(\star\otimes \succ,\succ\otimes\succ),
(\succ\otimes\circ,\succ\otimes\circ), \\
& (\prec\otimes\circ,\circ\otimes\succ),(\circ\otimes\prec,\circ\otimes\prec),
(\circ\otimes\circ,\circ\otimes\circ)\}. \notag
\end{eqnarray}
Note that the trialgebra contains an associative operation $\circ$
which is part of the splitting of the associative operation $\star$.
\smallskip
\\
4. (NS-algebra) The {\bf NS-algebra} of Leroux~\cite{Le2} is defined with
three binary operators
$\prec,\succ,\bullet$ that satisfy the relations
$$(x\prec y)\prec z=x\prec (y\star z),
(x\succ y)\prec z=x\succ (y\prec z),
(x\star y)\succ z=x\succ (y\succ z),$$
$$(x\star y)\bullet z+(x\bullet y)\prec z = x \succ (y\bullet z)
+x\bullet (y\star z)$$
for $x,y,z\in D$. Here $\star=\prec+\succ+\, \bullet$ gives an associative 
operation.
This is the \Loday algebra of type $(\dfgenb_N, \dfrelb_N)$
with $\dfgenb_N=\{\prec,\succ\,\bullet\}$ and
\begin{eqnarray}
\dfrelb_N=& \{(\prec\otimes\prec,\prec\otimes\star),
(\succ\otimes\prec,\succ\otimes\prec),
(\star\otimes \succ,\succ\otimes\succ), \mlabel{eq:NS}\\
& (\star\otimes\bullet+\bullet\otimes\prec,
\succ\otimes\bullet+\bullet\otimes \star)\}. \notag
\end{eqnarray}
5. ($L$-dipterous algebras)  An $L$-dipterous algebra~\cite{L-R3} is a $\bfk$-vector
space $A$ with two operations $\star$ and $\succ$ such that
\begin{equation}
\dfrelb_L=\{(\star\otimes \star, \star\otimes \star), (\star\otimes \succ,
\succ\otimes \succ), (\succ\otimes \star, \succ\otimes \star)\}.
\mlabel{eq:dip}
\end{equation}
An $L$-anti-dipterous algebra is a vector space $A$ with two operations
$\star$ and $\prec$ such that
\begin{equation}
\dfrelb_R=\{(\star\otimes \star, \star\otimes \star), (\prec\otimes\prec,
\prec\otimes\star), (\star\otimes \prec, \star\otimes \prec)\}.
\mlabel{eq:anti}
\end{equation}

\section{Products of \Loday operads}
\mlabel{sec:prod}
\subsection{Definitions}
Manin~\cite{Ma1,Ma2} defined two products for quadratic algebras,
called the {\bf black circle product} $\bullet$ and {\bf white
circle product} $\circ$. They were later generalized for operads~\cite{G-K}. 
We defined similar products for \Loday operads in terms of the types. 
The analogue of the black circle product was independently defined by 
Loday~\cite{Lo7} where he called it the {\bf square product} with 
the notation $\square$. We adopt his terminology and notation. We also call
the analogue of the white circle product the {\bf maltese product} with 
the notation $\maltese$. 
The justification of the names and notations will be given below. 
As the referee and Loday pointed out, the square product and maltese product 
are different from the black circle product and 
white circle product of Manin-Ginzburg-Kapranov. We will not go further in this 
direction since the later products will not be needed in this paper. A remark 
is given at the end of \S~\ref{sec:dual}. 

Let $(\dfgen_1, \dfrel_1)$ and
$(\dfgen_2,\dfrel_2)$ be the types of two \Loday operads. Define
$$ (\dfgen_1,\dfrel_1)\square (\dfgen_2,\dfrel_2)
=(\dfgen_1\otimes \dfgen_2, S_{(23)}(\dfrel_1\otimes \dfrel_2)),$$
$$ (\dfgen_1,\dfrel_1)\maltese (\dfgen_2,\dfrel_2)=
(\dfgen_1\otimes \dfgen_2, S_{(23)}(\dfrel_1\otimes \dfgen_2^{\otimes 2}+
    \dfgen_1^{\otimes 2} \otimes \dfrel_2)).$$
Here $S_{(23)}$ means the exchange of factor 2 and 3 in the tensor products. 
We see that relations of the square product (resp. maltese product) take the 
shape of a square (resp. a (Maltese) cross).
We will rephrase these products more precisely for later applications.

For $(\dfgen_1,\dfrel_1)$ and $(\dfgen_2,\dfrel_2)$ as above, and for
$\dfop^{(i)}\in \dfgen_i,\ i=1,2$, we use a column vector
$\dfpair{\dfop_1}{\dfop_2}$ to denote the tensor product
$\dfop_1\otimes \dfop_2\in \dfgen_1\otimes \dfgen_2$.
The purpose of this is to distinguish it from the tensor
product in $\dfrel_i\subseteq \dfgen_i^{\otimes 2}\oplus \dfgen_i^{\otimes 2}.$
Further, for
$f_i=(\dfoa_i\otimes\dfob_i, \dfoc_i\otimes\dfod_i)
\in \dfgen_i^{\otimes 2} \oplus \dfgen_i^{\otimes 2},\ i=1,2,$
define
$$\dfpair{f_1}{f_2}=
 \left (
\dfpair{\dfoa_1}{\dfoa_2}
\otimes\dfpair{\dfob_1}{\dfob_2},
\dfpair{\dfoc_1}{\dfoc_2}
\otimes\dfpair{\dfod_1}{\dfod_2}\right ) \in
(\dfgen_1\otimes \dfgen_2)^{\otimes 2}\oplus (\dfgen_1\otimes \dfgen_2)^{\otimes 2}.
$$
This extends by bilinearity to all
$f_i\in \dfgen_i^{\otimes 2} \oplus \dfgen_i^{\otimes 2},\ i=1,2.$
More precisely, elements of $\dfgen_i^{\otimes 2} \oplus \dfgen_i^{\otimes 2}$ are
finite sums of the form
$$f_1= \sum_j (\dfoa_{1,j}\otimes\dfob_{1,j}, \dfoc_{1,j}\otimes\dfod_{1,j}),\
\dfop_{1,j}^{(r)}\in \dfgen_1, \; r=1,\dots,4 $$ and
$$f_2= \sum_k (\dfoa_{2,k}\otimes\dfob_{2,k}, \dfoc_{2,k}\otimes\dfod_{2,k}),\
\dfop_{2,k}^{(r)}\in \dfgen_2, \; r=1,\dots,4.$$ We then define
$$\dfpair{f_1}{f_2}=
\sum_{j,k} \left ( \dfpair{\dfoa_{1,j}}{\dfoa_{2,k}}
\otimes\dfpair{\dfob_{1,j}}{\dfob_{2,k}},
\dfpair{\dfoc_{1,j}}{\dfoc_{2,k}}
\otimes\dfpair{\dfod_{1,j}}{\dfod_{2,k}}\right).
$$
We define two subspaces of $(\dfgen_1\otimes \dfgen_2)^{\otimes
2}\oplus (\dfgen_1\otimes \dfgen_2)^{\otimes 2}$ by
$$\dfrel_1\dftimes \dfrel_2=
\left \{ \dfpair{f_1}{f_2} \Big | f_i \in \dfrel_i,\ i=1,2 \right
\}, \quad \dfrel_1\dfotimes \dfrel_2= \left \{ \dfpair{f_1}{f_2}
\Big | f_1 \in \dfrel_1 {\rm\ or\ } f_2\in \dfrel_2 \right \}.
$$
So $\dfrel_1\dftimes \dfrel_2$ and $\dfrel_1\dfotimes \dfrel_2$
can be regarded as sets of relations for the
operator set
$\dfgen_1\otimes \dfgen_2$.

\begin{defn}
The {\bf square product} {\rm (}also called {\bf type product}{\rm )} 
of $(\dfgen_1,\dfrel_1)$ and $(\dfgen_2,\dfrel_2)$, 
denoted by $(\dfgen_1,\dfrel_1)\dftimes (\dfgen_2,\dfrel_2)$, 
is the type $(\dfgen_1\otimes \dfgen_2, \dfrel_1\dftimes \dfrel_2)$.
The {\bf maltese product} of $(\dfgen_1,\dfrel_1)$ and $(\dfgen_2,\dfrel_2)$, 
denoted by $(\dfgen_1,\dfrel_1)\dfotimes (\dfgen_2,\dfrel_2)$, 
is the type
$(\dfgen_1\otimes \dfgen_2, \dfrel_1\dfotimes \dfrel_2)$.
\end{defn}

Let $\dfgenb_i=\{\dfgene_{i,j}\}$ be a basis of $\dfgen_i$, $i=1,2$.
Then
$$ \dfgenb_1\dftimes \dfgenb_2:= \left\{ \dfpair{\dfgene_{1,j}}{\dfgene_{2,k}}
\right \}_{j,k}$$ is a basis of $\dfgen_1\otimes \dfgen_2$. Similarly,
let $\dfrelb_i=\{\dfrele_{i,j}\}$ be a basis of $\dfrel_i$,
$i=1,2$. Then
$$ \dfrelb_1\dftimes \dfrelb_2:= \left\{ \dfpair{\dfrele_{1,j}}{\dfrele_{2,k}}
\right \}_{j,k}$$
is a basis of $\dfrel_1\dftimes \dfrel_2$ and takes the shape of a box (matrix).
Let $\vec{f}_i$ be a basis of $\dfgen_i^{\otimes 2}\oplus \dfgen_i^{\ot 2}$, $i=1,2$.
Then a spanning set  of $\dfrel_1 \dfotimes \dfrel_2$ is given by
$$ (\dfrelb_1\dftimes \vec{f}_2) \bigcup\, (\vec{f}_1 \dftimes \dfrelb_2). $$
This set is not linearly independent. The union even has overlap
when $\vec{f}_i$ is extended from $\dfrelb_i$. Then the union takes the shape 
of a cross. We will mostly consider the type product $\dftimes$. 
Concerning the duality (or the lack of it) between the products $\dftimes$ and 
$\dfotimes$, see \S\ref{sec:dual} for details.

\subsection{Basic properties}
\subsubsection{Transposes and isomorphisms}
Let $(\dfgen_1\otimes \dfgen_2,\dfrel_1\dftimes \dfrel_2)$ be
the type product of two \Loday algebras.
We define the {\bf transpose} of $(\dfgen_1\otimes\dfgen_2,\dfrel_1\dftimes\dfrel_2)$
by $(\dfgen_2\otimes\dfgen_1,\dfrel_2\dftimes\dfrel_1)$.
It is obtained by changing the two factors in
$\dfpair{\dfop_1}{\dfop_2}$ throughout
$\dfgen_1\otimes\dfgen_2$ and
$\dfrel_1\dftimes\dfrel_2$.

\begin{lemma}
Let $(\dfgen_i,\dfrel_i)$, $i=1,2$, be the types of 
two \Loday operads with associative operations $\star_i$.
\begin{enumerate}
\item
$(\dfgen_1\otimes \dfgen_2,\dfrel_1\dftimes \dfrel_2)$ is an \Loday operad with 
associative operation $\dfpair{\star_1}{\star_2}$.
\item $(\dfgen_1\ot \dfgen_2, \dfrel_1\dftimes \dfrel_2)$ is isomorphic to its
{ transpose} $(\dfgen_2\ot \dfgen_1, \dfrel_2\dftimes \dfrel_1)$.
\item
If $(\dfgen_1,\dfrel_1)$ is isomorphic to $(\dfgen'_1,\dfrel'_1)$
and $(\dfgen_2,\dfrel_2)$ is isomorphic to
$(\dfgen'_2,\dfrel'_2)$,
then $(\dfgen_1\otimes \dfgen_2,\dfrel_1\dftimes \dfrel_2)$ is
isomorphic to
$(\dfgen'_1\otimes \dfgen'_2,\dfrel'_1\dftimes \dfrel'_2)$.
\end{enumerate}
\mlabel{lem:proj}
\end{lemma}
\begin{proof}
(1). Since $\star_i$ is associative, we have
$(\star_i\otimes \star_i, \star_i\ot \star_i)\in \dfrel_i$.
Thus the product
$\dfpair{(\star_1\otimes \star_1, \star_1\ot \star_1)}
{(\star_2\otimes \star_2, \star_2\ot \star_2)}$ is in $\dfrel_1\dftimes \dfrel_2$.
But this product is just
$\left (\dfpair{\star_1}{\star_2}\ot \dfpair{\star_1}{\star_2},
\dfpair{\star_1}{\star_2}\ot \dfpair{\star_1}{\star_2} \right )$.
So $\dfpair{\star_1}{\star_2}$ is associative.

(2). The linear map $\dfgen_1\otimes \dfgen_2\to \dfgen_2\otimes \dfgen_1$
sending $\dfpair{\dfgene_1}{\dfgene_2}$ to $\dfpair{\dfgene_2}{\dfgene_1}$
is bijective that sends $\dfpair{\star_1}{\star_2}$ to $\dfpair{\star_2}{\star_1}$
and induces a bijective linear map
$\dfrel_1\otimes \dfrel_2\to \dfrel_2\otimes \dfrel_1$. This proves the claim.

The proof of (3) is similar. \end{proof}

\subsubsection{Tensor products}

\begin{prop}
Let $(D_1,\dfgen_1,\dfrel_1)$ and $(D_2,\dfgen_2,\dfrel_2)$
be \Loday algebras of type $(\dfgen_1,\dfrel_1)$
and $(\dfgen_2,\dfrel_2)$ respectively.
For $a_1\otimes a_2, b_1\otimes b_2\in D_1\otimes D_2$
and $\dfop_1\in \dfgen_1,\dfop_2 \in \dfgen_2$,
define
$$(a_1\otimes a_2) \dfpair{\dfop_1}{\dfop_2} (b_1\otimes b_2)
=(a_1\dfop_1 b_1) \otimes (a_2 \dfop_2 b_2).$$
This defines an \Loday algebra of type
$(\dfgen_1 \otimes \dfgen_2, \dfrel_1\dftimes \dfrel_2)$
on $D_1\otimes D_2$.
\end{prop}
For example, when both $(\dfgen_1,\dfrel_1)$ and
$(\dfgen_2,\dfrel_2)$ are of the type of a dendriform dialgebra,
then $D_1\otimes D_2$ is a quadri-algebra. See~\cite[1.5]{A-L}.
When both types are of a dendriform trialgebra, then $D_1\otimes
D_2$ is an ennea-algebra. When the two types are of trialgebra and
of NS-algebra, respectively, then $D_1\otimes D_2$ is a
dendriform-Nijenhuis algebra~\cite{Le2}.

\begin{proof}
We only need to verify that the relations in
$\dfrel_1\dftimes \dfrel_2$ are satisfied by $D_1\otimes D_2$.

Recall that $\dfrel_1\dftimes \dfrel_2$ consists of elements
of the form $\dfpair{f_1}{f_2}$
with $f_i\in \dfrel_i,\ i=1,2$.
Further,  if
$$f_i=\sum_{j_i} \left (
\dfoa_{i,j_i}\otimes\dfob_{i,j_i},
\dfoc_{i,j_i}\otimes\dfod_{i,j_i}
\right )\in \dfrel_i,\ i=1,2,$$
then
$$
\dfpair{f_1}{f_2}=\sum_{j_1,j_2} \left (
\dfpair{\dfoa_{1,j_1}}{\dfoa_{2,j_2}}
\otimes\dfpair{\dfob_{1,j_1}}{\dfob_{2,j_2}},
\dfpair{\dfoc_{1,j_1}}{\dfoc_{2,j_2}}
\otimes\dfpair{\dfod_{1,j_1}}{\dfod_{2,j_2}}\right ).
$$
For $x_i,y_i,z_i\in D_i,\ i=1,2$, we have
\allowdisplaybreaks{
\begin{eqnarray*}
&&\sum_{j_1,j_2}
\left ((x_1\otimes x_2) \dfpair{\dfoa_{1,j_1}}{\dfoa_{2,j_2}}
(y_1\otimes y_2)\right )
\dfpair{\dfob_{1,j_1}}{\dfob_{2,j_2}} (z_1\otimes z_2)\\
&=& \sum_{j_1,j_2}
\left ((x_1\dfoa_{1,j_1} y_1)\dfob_{1,j_1} z_1\right ) \otimes
\left ((x_2 \dfoa_{2,j_2} y_2) \dfob_{2,j_2} z_2 \right )
\\
&=&
 \left (\sum_{j_1}
 (x_1\dfoa_{1,j_1} y_1)\dfob_{1,j_1} z_1 \right )\otimes
\left ( \sum_{j_2}
 (x_2 \dfoa_{2,j_2} y_2) \dfob_{2,j_2} z_2 \right )\\
&=&
\left (\sum_{j_1}
 x_1\dfoc_{1,j_1} (y_1\dfod_{1,j_1} z_1) \right )\otimes
\left ( \sum_{j_2}
 x_2 \dfoc_{2,j_2} (y_2 \dfod_{2,j_2} z_2) \right )\\
&=&
\sum_{j_1,j_2}
 x \dfpair{\dfoc_{1,j_1}}{\dfoc_{2,j_2}}
\left ( y\dfpair{\dfod_{1,j_1}}{\dfod_{2,j_2}} z
\right ).
\end{eqnarray*}
}
This is what we want.
\end{proof}

\subsection{Examples}
We now show that some recent generalizations of dendriform dialgebras are products of
more basic algebras. Therefore, their generators and relations can be easily described.

\subsubsection{Quadri-algebra}
The quadri-algebra of Aguiar and Loday~\cite{A-L} is defined by four
binary operations $\{\nearrow\,,\nwarrow\,,\searrow\,,\swarrow\,\}$ and
9 relations.
Using the auxiliary operations
$$\wedge=\,\nearrow+\nwarrow\,, \vee=\,\searrow+\swarrow\,, \prec\,=\,\nwarrow+\swarrow\,,
\succ\,=\,\nearrow+\searrow\,, \star = \wedge+\vee=\,\prec+\succ\,,$$
the 9 relations are given in the following $3\times 3$ matrix.

\begin{equation}
\raisebox{25pt}{\xymatrix@C=5pt@R=5pt{
{\mbox{\footnotesize $(x\nw y)\nw z =x\nw(y\star z)$}} & &
{\mbox{\footnotesize $(x\ne y)\nw z = x\ne(y\west z)$}} &  &
{\mbox{\footnotesize $(x\north y)\ne z=x\ne(y\east z)$}}  \\
 {\mbox{\footnotesize $(x\sw y)\nw z =x\sw(y\north z)$}} & &
 {\mbox{\footnotesize $(x\se y)\nw z = x\se(y\nw z)$}} &  &
{\mbox{\footnotesize $(x\south y)\ne z=x\se(y\ne z)$}}             \\
 {\mbox{\footnotesize $(x\west y)\sw z=x\sw(y\south z)$}} & &
 {\mbox{\footnotesize $(x\east y)\sw z = x\se(y\sw z)$}} &  &
{\mbox{\footnotesize $(x\star y)\se z =x\se(y\se z)$}}
}}
\mlabel{eq:quadrie}
\end{equation}
\smallskip

We will put quadri-algebra in the context of type products.
Recall that the dendriform dialgebra of Loday is of type $(\dfgenb_D,\dfrelb_D)$
with $\dfgenb_D=\{\prec,\succ\}$ and
$$\dfrelb_D=\{(\prec\otimes\prec,\prec\otimes\dft),
(\succ\otimes\prec,\succ\otimes\prec),
(\dft\otimes\succ,\succ\otimes\succ)\}$$
where $\dft=\, \succ+\prec$.

\begin{prop}
The quadri-algebra is isomorphic to the \Loday operad of type
$$(\dfgen_D,\dfrel_D)
\dftimes (\dfgen_D,\dfrel_D)=(\dfgen_D\otimes\dfgen_D,
\dfrel_D\dftimes\dfrel_D).$$
\mlabel{pp:quadri}
\end{prop}
Since the two factors of the product are the same, the transpose
of a quadri-algebra is the same algebra. Since the opposite of a
dialgebra is still a dialgebra, by Lemma~\ref{lem:proj}, if we
exchange one or both pairs of $\prec$ and $\succ$, we still obtain
a quadri-algebra. The opposite quadri-algebra in~\cite{A-L} is
obtained when both pairs are exchanged. As a consequence, applying
every element from the dihedral group $D_4$ of order 8 to a
quadri-algebra again gives a quadri-algebra.

\begin{proof}
We first define a bijection between the binary operations of the quadri-algebra and
the binary operations
$$\dfgenb_D\otimes \dfgenb_D=\left \{ \dfpair{\prec}{\prec},\dfpair{\prec}{\succ},
\dfpair{\succ}{\prec}, \dfpair{\succ}{\succ} \right \}$$ of the
type product which is given in the following table.

\begin{equation}
\begin{array}{ccc}
\nwarrow\leftrightarrow\dfpair{\prec}{\prec}, & \nearrow\leftrightarrow\dfpair{\prec}{\succ}, &
\wedge\leftrightarrow\dfpair{\prec}{\dft} \vspace{.3cm}\\
\swarrow\leftrightarrow\dfpair{\succ}{\prec}, & \searrow\leftrightarrow\dfpair{\succ}{\succ}, &
\vee\leftrightarrow\dfpair{\succ}{\dft} \vspace{.3cm} \\
\prec\leftrightarrow\dfpair{\dft}{\prec}, & \succ\leftrightarrow\dfpair{\dft}{\succ}, &
\star\leftrightarrow\dfpair{\dft}{\dft}
\end{array}
\mlabel{eq:quadri}
\end{equation}
Here the entries
on the third row and column are defined to be the sums along the corresponding projections.

To help visualizing this bijection, the reader can imagine the
$xy$ plane as sitting in the 3-dimensional coordinate system in the usual way.
Thus the plane is laying flat with the $x$-axis pointing outwards to the
reader and the $y$-axis pointing to the right. Also $\prec$ (resp.
$\succ$) points to the negative (resp. positive) direction of the axis.
Then, for example, the northwest arrow $\nwarrow$ should be visualized not
pointing up and left, but rather inward and left to the third
quadrant --- the negative
direction in both the $x$ and $y$ coordinates. This agrees with
the meaning of $\dfpair{\prec}{\prec}$.

Under the bijection in Eq.~(\ref{eq:quadri}), the relation matrix~(\ref{eq:quadrie})
is sent to
\begin{eqnarray*}
&&\hspace{-1.3cm}
\begin{array}{ccc}
(x\dfrr y)\dfrr z=x\dfrr(y\dftt z), &
(x\dfrl y)\dfrr z=x\dfrl (y\dftr z), &
(x\dfrt y)\dfrl z=x\dfrl (y\dftl z),\vspace{.3cm}\\
(x\dflr y)\dfrr z=x\dflr(y\dfrt z), &
(x\dfll y)\dfrr z=x\dfll(y\dfrr z), &
(x\dflt y)\dfrl z=x\dfll(y\dfrl z),
\vspace{.3cm}\\
(x\dftr y)\dflr z=x\dflr(y\dflt z), &
(x\dftl y)\dflr z=x\dfll(y\dflr z), &
(x\dftt y)\dfll z=x\dfll(y\dfll z)
\vspace{.3cm}
\end{array}
\end{eqnarray*}
which is simply the matrix of $\dfrelb_D \dftimes \dfrelb_D$.
\end{proof}

\subsubsection{Ennea-algebra}
The ennea-algebra (or $1$-ennea-algebra) of Leroux~\cite{Le1} has 9 binary operations
$$
\nwarrow\,, \uparrow\,,\nearrow\,, \prec\,,\circ\,,\succ\,,\swarrow\,,
\downarrow\,, \searrow$$
and 49 relations.

\begin{prop}
The ennea-algebra is isomorphic to the type product
$$(\dfgen_T,\dfrel_T)\dftimes (\dfgen_T,\dfrel_T)$$
where $(\dfgen_T,\dfrel_T)$ is the type of the dendriform
trialgebra of Loday and Ronco. \mlabel{pp:ennea}
\end{prop}
Again it is apparent that exchanging $\dfpair{\dfop_1}{\dfop_2}
\leftrightarrow \dfpair{\dfop_2}{\dfop_1}$ gives the transpose of
the relation matrix. So the transpose of an ennea-algebra is the
same algebra. It is also clear that the opposite algebra is the
self-product of the opposite trialgebra, and elements in $D_4$
give algebras isomorphic to the ennea-algebra.

\begin{proof}
Recall that the dendriform trialgebra of Loday and Ronco~\cite{L-R2}
is of type $(\dfgenb_T,\dfrelb_T)$ with
$\dfgenb_T=\{\prec,\succ\,\circ\}$ and
\begin{eqnarray*}
\dfrelb_T=& \{(\prec\otimes\prec,\prec\otimes\star),
(\succ\otimes\prec,\succ\otimes\prec),
(\star\otimes \succ,\succ\otimes\succ),
(\succ\otimes\circ,\succ\otimes\circ), \\
& (\prec\otimes\circ,\circ\otimes\succ),(\circ\otimes\prec,\circ\otimes\prec),
(\circ\otimes\circ,\circ\otimes\circ)\}.
\end{eqnarray*}

Similar to the quadri-algebra, we first give a bijection
between the operations of Leroux and those in $\dfgen_T\otimes \dfgen_T$
in the following table. The entries
on the fourth row and column are defined to be the sums along the corresponding projections.

$$\begin{array}{cccc}
\nwarrow\leftrightarrow\dfpair{\prec}{\prec}, & \uparrow \leftrightarrow \dfpair{\prec}{\circ}, &
\nearrow\leftrightarrow\dfpair{\prec}{\succ}, &
\wedge\leftrightarrow\dfpair{\prec}{\dft} \vspace{.3cm}\\
\prec\leftrightarrow\dfpair{\circ}{\prec}, & \circ\leftrightarrow\dfpair{\circ}{\circ},&
\succ\leftrightarrow\dfpair{\circ}{\succ}, & \star\leftrightarrow\dfpair{\circ}{\dft} \vspace{.3cm}\\
\swarrow\leftrightarrow\dfpair{\succ}{\prec}, & \downarrow\leftrightarrow\dfpair{\succ}{\circ}, &
\searrow\leftrightarrow\dfpair{\succ}{\succ}, &
\vee\leftrightarrow\dfpair{\succ}{\dft} \vspace{.3cm} \\
\triangleleft\leftrightarrow\dfpair{\dft}{\prec}, & \bar{\circ}\leftrightarrow\dfpair{\dft}{\circ},
& \triangleright\leftrightarrow\dfpair{\dft}{\succ}, &
\bar{\star}\leftrightarrow\dfpair{\dft}{\dft}
\end{array}
$$
Applying this bijection to the 49 relations for the ennea-algebra in \cite{Le1},
we see that the relations become the entries of the $7\times 7$ matrix
$ \dfrelb_T \dftimes \dfrelb_T.$
\end{proof}

\subsubsection{Dendriform-Nijenhuis algebra}

The dendriform-Nijenhuis algebra~\cite{Le2}
is equipped with 9 binary operations
\begin{equation}
\nearrow\,,\ \searrow\,,\ \swarrow\,,\ \nwarrow\,,\ \uparrow\,,\ \downarrow\,,\
\tilde{\prec}\,, \ \tilde{\succ}\,,\ \tilde{\bullet}
\mlabel{eq:dN}
\end{equation}
satisfying 28 relations.

\begin{prop}
The dendriform-Nijenhuis algebra is isomorphic to the \Loday algebra of type
$(\dfgen_T,\dfrel_T)\dftimes (\dfgen_N,\dfrel_N)$ where
$(\dfgen_T,\dfrel_T)$ is the type of dendriform trialgebra and
$(\dfgen_N,\dfrel_N)$ is the type of a NS-algebra in
Ep.~(\ref{eq:NS}).
\end{prop}
\begin{proof}
We just give a correspondence  between the binary operation of
Leroux and ours. To distinguish operators in the dendriform
trialgebra and NS-algebra, we denote $\dfgen_N=\{<,>,\bullet\}$
and $\diamondsuit=<+>+\,\bullet$. A bijection between the binary
operations in Eq.~(\ref{eq:dN}) and $\dfgenb_T\otimes \dfgenb_N$ is
given by
$$\begin{array}{cccc}
\nwarrow\leftrightarrow\dfpair{\prec}{<}, & \uparrow \leftrightarrow \dfpair{\prec}{\bullet}, &
\nearrow\leftrightarrow\dfpair{\prec}{>}, &
\wedge\leftrightarrow\dfpair{\prec}{\diamondsuit} \vspace{.3cm}\\
\tilde{\prec}\leftrightarrow\dfpair{\circ}{<}, & \tilde{\bullet}\leftrightarrow\dfpair{\bullet}{\circ},&
\tilde{\succ}\leftrightarrow\dfpair{\circ}{>}, & \tilde{\star}\leftrightarrow\dfpair{\circ}{\diamondsuit} \vspace{.3cm}\\
\swarrow\leftrightarrow\dfpair{\succ}{<}, & \downarrow\leftrightarrow\dfpair{\succ}{\bullet}, &
\searrow\leftrightarrow\dfpair{\succ}{>}, &
\vee\leftrightarrow\dfpair{\succ}{\diamondsuit} \vspace{.3cm} \\
\triangleleft\leftrightarrow\dfpair{\dft}{<}, & \bar{\bullet}\leftrightarrow\dfpair{\dft}{\bullet},
& \triangleright\leftrightarrow\dfpair{\dft}{>}, &
\bar{\star}\leftrightarrow\dfpair{\dft}{\diamondsuit}
\end{array}
$$
The operations on the fourth row and column are defined to be the
sum along the corresponding projections. Then the relations in a
dendriform-Nijenhuis algebra~\cite{Le2} is identified with the
$7\times 4$ matrix $\dfrelb_T\dftimes \dfrelb_N$.
\end{proof}

\subsubsection{Octo-algebra}
The octo-algebra of Leroux~\cite{Le3} is isomorphic to the product
$(\dfgen_Q\otimes \dfgen_D,\dfrel_Q\dftimes \dfrel_D)$ where
$(\dfgen_Q,\dfrel_Q)$ is the type of quadri-algebra and
$(\dfgen_D,\dfrel_D)$ is the type of dendriform dialgebra. As we
will see later in \S\ref{ss:power}, it is also the third power of
dendriform dialgebra defined there. We will give details there
(Proposition~\ref{pp:octo}).

\subsubsection{Type $M_1$ and $M_2$ algebras}
Type $M_1$ and $M_2$ algebras were introduced in~\cite{Le3} as
expansions of dipterous and anti-dipterous algebras.
Let $(\dfgen_i,\dfrel_i), i=1,2,$ both be the type of the $L$-dipterous
algebra. So we have
$\dfgen_i=\{\star_i,\succ_i\}$ with
$$\dfrel_i=\{(\star_i\otimes \star_i, \star_i\otimes \star_i),
(\star_i\otimes \succ_i, \succ_i\otimes \succ_i), (\succ_i\otimes
\star_i, \succ_i\otimes \star_i)\},\; i=1,2.$$ Then the $M_2$
algebra is isomorphic to the product $(\dfgen_1\otimes
\dfgen_2,\dfrel_1\dftimes \dfrel_2)$. The correspondence between
the four binary operations of $M_2$ and $\dfgenb_1\otimes
\dfgenb_2$ is
$$\bullet_1=\dfpair{\succ}{\succ}, \bullet_2=\dfpair{\succ}{\star},
\bullet_3=\dfpair{\star}{\succ}, \bullet_4=\dfpair{\star}{\star}.$$
We skip the subscripts $i=1,2,$ since it is clear from the context.

Similarly, The $M_1$ algebra in~\cite{Le3} is the product of
the $L$-anti-dipterous algebra and $L$-dipterous algebra.

\subsection{Powers of an \Loday operad}
\mlabel{ss:power}
Inductively, we define the  type product of any finite number of
types of \Loday operads: given $(\dfgen_i,\dfrel_i), 1\leq i\leq n$,
define
$$\dfTimes_{i=1}^n \left (\dfgen_i,\dfrel_i\right):=
\left (\dfTimes_{i=1}^{n-1} \left (\dfgen_i,\dfrel_i
\right )\right )\dftimes \left (\dfgen_n,\dfrel_n\right).$$
In particular, we define
the powers of a specific type of \Loday operad.
For simplicity, we only describe the powers of a dendriform trialgebra.
\begin{defn}
A $N$-th power of the trialgebra is a $\bfk$-vector space $D$
equipped with $3^N$ binary operations
$$\dftri{\dfop_1}{\cdots}{\dfop_N},\
\dfop_i\in\{\dfl_i,\dfr_i,\dfc_i\},\ 1\leq i\leq N$$
such that, for any choice of
$(\dfoa_i\otimes\dfob_i,\dfoc_i\otimes\dfod_i)$ in
$$ \dfrel_i:=
\left \{\begin{array}{l}
(\dfr_i\otimes\dfr_i,\dfr_i\otimes\dft_i),(\dfl_i\otimes\dfr_i,
\dfl_i\otimes\dfr_i),
(\dft_i\otimes\dfl_i,\dfl_i\otimes\dfl_i),(\dfl_i\otimes\dfc_i,
\dfl_i\otimes\dfc_i),\\
(\dfr_i\otimes\dfc_i,\dfc_i\otimes\dfl_i),(\dfc_i\otimes\dfr_i,
\dfc_i\otimes\dfr_i),
(\dfc_i\otimes\dfc_i,\dfc_i\otimes\dfc_i)
\end{array} \right \}$$
as $1\leq i\leq N$,
we have
$$\left (x \dftri{\dfoa_1}{\cdots}{\dfoa_N} y \right )\dftri{\dfob_1}{\cdots}{\dfob_N}z
=x \dftri{\dfoc_1}{\cdots}{\dfoc_N} \left (y \dftri{\dfod_1}{\cdots}{\dfod_N}z \right).
$$
Here $\dft_i=\dfl_i+\dfr_i+\dfc_i$. When $\dfc_i=0$ for $1\leq
i\leq N$, we call the algebra the $N$-th power of the dendriform
dialgebra.
\end{defn}

Of course, the dendriform dialgebra and quadri-algebra (resp.
dendriform trialgebra and ennea-algebra) are just the first and
second power of dendriform dialgebra (resp. trialgebra). The
octo-algebra introduced by Leroux~\cite{Le3} is the third power of
the dendriform dialgebra (see below).

\subsection{Examples}
We now give examples of algebras with powers greater than two and
more generally, of products with more than two factors.

\subsubsection{Octo-algebras}
The octo-algebra was introduced by Leroux~\cite{Le3}.
It is defined using 8 operations
$$\nearrow_i,\nwarrow_i,\swarrow_i,\searrow_{\,i}\,,\ i=1,2$$
and 27 relations.

\begin{prop}
The octo-algebra is the third power of the dendriform dialgebra
$(\dfgenb_D,\dfrelb_D)$.
\mlabel{pp:octo}
\end{prop}
\begin{proof}
The following table gives a correspondence between these 8 operations and
the operations in $\dfgenb_D \otimes \dfgenb_D \otimes \dfgenb_D
=\dfgenb_Q \otimes \dfgenb_D$ with $\dfgenb_Q$ the generators of the quadri-algebra.
As in the case of the quadri-algebra, each operation on the third row and third column
in the first two blocks is the sum of the corresponding projection, and each
operation in the third block is the sum of the corresponding operations on the
first and second blocks.
$$\begin{array}{ccc}
\nwarrow_1\leftrightarrow\dftri{\prec}{\prec}{\scriptscriptstyle{\prec}},
&
\nearrow_1\leftrightarrow\dftri{\prec}{\succ}{\scriptscriptstyle{\prec}},
&
\wedge_1\leftrightarrow\dftri{\scriptscriptstyle{\prec}}{\dft}{\scriptscriptstyle{\prec}}, \vspace{.3cm}\\
\swarrow_1\leftrightarrow\dftri{\succ}{\prec}{\scriptscriptstyle{\prec}},
&
\searrow_1\leftrightarrow\dftri{\succ}{\succ}{\scriptscriptstyle{\prec}},
&
\vee_1\leftrightarrow\dftri{\succ}{\dft}{\scriptscriptstyle{\prec}}, \vspace{.3cm} \\
\prec_1\leftrightarrow\dftri{\dft}{\prec}{\scriptscriptstyle{\prec}},
&
\succ_1\leftrightarrow\dftri{\dft}{\succ}{\scriptscriptstyle{\prec}},
&
\star_1\leftrightarrow\dftri{\dft}{\dft}{\scriptscriptstyle{\prec}},
\vspace{.3cm}\\
\hline
\end{array}
$$
$$\begin{array}{ccc}
\nwarrow_2\leftrightarrow\dftri{\prec}{\prec}{\scriptscriptstyle{\succ}},
&
\nearrow_2\leftrightarrow\dftri{\prec}{\succ}{\scriptscriptstyle{\succ}},
&
\wedge_2\leftrightarrow\dftri{\prec}{\dft}{\scriptscriptstyle{\succ}}, \vspace{.3cm}\\
\swarrow_2\leftrightarrow\dftri{\succ}{\prec}{\scriptscriptstyle{\succ}},
&
\searrow_2\leftrightarrow\dftri{\succ}{\succ}{\scriptscriptstyle{\succ}},
&
\vee_2\leftrightarrow\dftri{\succ}{\dft}{\scriptscriptstyle{\succ}}, \vspace{.3cm} \\
\prec_2\leftrightarrow\dftri{\dft}{\prec}{\scriptscriptstyle{\succ}},
&
\succ_2\leftrightarrow\dftri{\dft}{\succ}{\scriptscriptstyle{\succ}},
&
\star_2\leftrightarrow\dftri{\dft}{\dft}{\scriptscriptstyle{\succ}},
\vspace{.3cm}\\
\hline
\end{array}
$$
$$\begin{array}{ccc}
\nwarrow_{12}\leftrightarrow\dftri{\prec}{\prec}{\scriptscriptstyle{\dft}},
&
\nearrow_{12}\leftrightarrow\dftri{\prec}{\succ}{\scriptscriptstyle{\dft}},
&
\bigwedge_{12}\leftrightarrow\dftri{\prec}{\dft}{\scriptscriptstyle{\dft}}, \vspace{.3cm}\\
\swarrow_{12}\leftrightarrow\dftri{\succ}{\prec}{\scriptscriptstyle{\dft}},
&
\searrow_{12}\leftrightarrow\dftri{\succ}{\succ}{\scriptscriptstyle{\dft}},
&
\bigvee_{12}\leftrightarrow\dftri{\succ}{\dft}{\scriptscriptstyle{\dft}}, \vspace{.3cm} \\
\ll\leftrightarrow\dftri{\dft}{\prec}{\scriptscriptstyle{\dft}}, &
\gg\leftrightarrow\dftri{\dft}{\succ}{\scriptscriptstyle{\dft}}, &
\bar{\star}\leftrightarrow\dftri{\dft}{\dft}{\scriptscriptstyle{\dft}},
\vspace{.3cm}
\end{array}
$$

Then the 27 axioms of an octo-algebra in~\cite{Le3} correspond to entries in
$$\dfrelb_D\dftimes \dfrelb_D \dftimes \dfrelb_D
= \dfrelb_Q \dftimes \dfrelb_D.$$
Here $\dfrelb_Q$ is the relation vector of the quadri-algebra.
\end{proof}

Recall that the opposite of the dialgebra is still a dialgebra.
Considering also permutations of the three coordinates in
$\dftri{\dfop_1}{\dfop_2}{\scriptscriptstyle{\dfop_3}}$, we find
that each element in the group of rigid motions of the cube gives
an algebra isomorphic to the octo-algebra.

\subsubsection{The di-dipterous-anti-dipterous algebra}
It is easy to get new structures. For example, taking the product
of a dendriform dialgebra, a $L$-dipterous algebra and a
$L$-anti-dipterous algebra gives the di-dipterous-anti-dipterous
algebra with 8 operations $\dfgenb_D\ot \dfgenb_L \ot \dfgenb_R$
and 27 relations $\dfrelb_D\ot \dfrelb_L \ot \dfrelb_R$.


\section{Duality of \Loday operads}
\mlabel{sec:dual}
We demonstrate how the description of operads in terms of their types can be 
used to give their duals.  We will also give some examples. 

\subsection{Definitions}
For an \Loday operad $\calp=(\dfgen,\dfrel)$, the dual operad is defined as follows. 
See~\cite[B.2]{Lo4} for further details. 

Let $\check{\dfgen}:=\Hom(\dfgen,\bfk)$ be the dual space of $\dfgen$, giving the natural 
pairing 
$$ \langle\ ,\  \rangle_\dfgen: \dfgen \times \check{\dfgen} \to \bfk.$$ 
Then via
$$\Hom(\dfgen^{\ot 2},\bfk) \cong \Hom(\dfgen,\check{\dfgen}) \cong 
\Hom(\dfgen,\bfk)\ot \check{\dfgen} \cong \check{\dfgen}^{\ot 2},$$ 
we get a natural (perfect) pairing 
\begin{equation}
\langle\ ,\  \rangle_{\Omega^{\ot 2}}: \dfgen^{\otimes 2} \times \check{\dfgen}^{\otimes 2} \to \bfk, 
\langle x\otimes y, a\otimes b\rangle_{\Omega^{\otimes 2}} 
=\langle x, a \rangle_\dfgen\ \langle y, b\rangle_\dfgen.
\mlabel{eq:pair1}
\end{equation}
Further the isomorphism 
$$ \Hom(\Omega^{\ot 2} \oplus \Omega^{\ot 2}, \bfk)\cong 
    \Hom(\Omega^{\ot 2},\bfk)\oplus \Hom(\Omega^{\ot 2},\bfk) 
    \cong \check{\Omega}^{\ot 2} \oplus \check{\Omega}^{\ot 2}$$
gives a perfect pairing 
\begin{equation}
\langle\ ,\  \rangle_{2\Omega^{\ot 2}}: 
(\dfgen^{\otimes 2} \oplus \dfgen^{\otimes 2}) \times 
    (\check{\dfgen}^{\otimes 2} \oplus \check{\dfgen}^{\otimes 2}) \to \bfk 
\mlabel{eq:pair2}
\end{equation}
by 
$$\langle (\alpha,\beta), (\gamma,\delta)\rangle_{2\Omega^{\ot 2}} 
= \langle \alpha,\gamma\rangle_{\Omega^{\otimes 2}}
-\langle \beta, \delta\rangle_{\Omega^{\otimes 2}}, \alpha,\beta\in \dfgen^{\otimes 2},\ 
\gamma,\delta\in\check{\dfgen}^{\otimes 2}.$$
More precisely, let $\{x_i\}$ be a basis of $\dfgen$ with dual basis $\{\check{x}_i\}$ 
in $\check{\dfgen}$. Then 
$$\langle (x_i\otimes x_j, x_k\otimes x_\ell), (\check{x}_s\otimes \check{x}_t, 
\check{x}_u\otimes \check{x}_v) \rangle_{2\Omega^{\ot 2}} = 
\delta_{i,s}\delta_{j,t}-\delta_{k,u}\delta_{\ell,v}.$$

We now define $\dfrel^{\bot}$ to be the annihilator of 
$\dfrel\subseteq \dfgen^{\otimes 2} \oplus \dfgen^{\otimes 2}$ in $\check{\dfgen}^{\otimes 2}\oplus 
\check{\dfgen}^{\otimes 2}$ under the pairing $\langle\ , \rangle_{2\Omega^{\ot 2}}$. 
We call $\calp^{\,!}:=(\check{\dfgen},\dfrel^{\bot})$ the {\bf dual operad} of 
$\calp=(\dfgen,\dfrel)$ which is the Koszul dual in our special case. 
It follows from the definition that $(\calp^{\,!})^!=\calp$.

\subsection{Examples}
\subsubsection{Dual operad of the dendriform dialgebra}
\mlabel{sec:dd}
This duality is given by Loday~\cite[Proposition 8.3]{Lo4}. 
Let $(\dfgen_D, \dfrel_D)$ be the type of the operad for the dendriform dialgebra. 
Let $\{\dashv,\vdash\}\in \check{\dfgen}_D$ be the dual basis of $\{\prec, \succ\}$ 
(in this order). Then $\dfrel_D^\bot$ is given by 
\begin{equation} \{ (\dashv\otimes \dashv, \dashv\otimes \dashv), 
(\vdash\otimes \vdash, \vdash\otimes \vdash), 
(\dashv \otimes\dashv, \dashv \otimes \vdash), 
(\vdash \otimes \dashv, \vdash \otimes \dashv),
(\dashv \otimes \vdash, \vdash \otimes \vdash)\}.
\mlabel{eq:asdi}
\end{equation}
It is called the associative dialgebra $(\dfgen_{AD},\dfrel_{AD})$. 
There are many associative binary operations in $\dfrel_{AD}$, such as  
$\vdash$, $\dashv$ or their linear combinations. 

\subsubsection{Dual operad of the dendriform trialgebra}
It is proved by Loday and Ronco~\cite[Theorem 3.1]{L-R2} that the dual operad of 
the dendriform trialgebra $(\dfgen_T,\dfrel_T)$ is the associative trialgebra with 
three generators and 11 relations. 

\subsubsection{Dual operad of the NS operad} 
We now find the dual of the NS-algebra. 
The NS-algebra of Leroux~\cite{Le2} is defined with 
three binary operators 
$\dfgenb_N=\{\prec,\succ,\bullet\}$ that satisfy the relations  
$$\dfrelb_N=\big \{(x\prec y)\prec z\oeq{1} x\prec (y\star z), 
(x\succ y)\prec z\oeq{2} x\succ (y\prec z), $$
$$(x\star y)\succ z \oeq{3} x\succ (y\succ z),
(x\star y)\bullet z+(x\bullet y)\prec z \oeq{4} x \succ (y\bullet z) 
+x\bullet (y\star z) \big \}$$
for $x,y,z\in D$. The labels on the equations are for later reference. 
 Here $\star=\prec+\succ+\bullet$. 
Let $\{\dashv,\vdash,\circ\}$ be the basis of $\check{\dfgen}_N$ that is dual 
to $\dfgenb_N$. 
We note that each term from $\dfgenb_N^{\otimes 2}$ occurs exactly once on 
the left hand side of the equations in $\dfrelb_N$ and exactly once on the right 
hand side. Thus by inspection, we found that the following elements in 
$\check{\dfgen}_N^{\otimes 2} \oplus \check{\dfgen}_N^{\otimes 2}$ are perpendicular 
to $\dfrelb_N$ and thus are in $\dfrel_N^\bot$. 
\begin{eqnarray*}
&&\hspace{-.7cm} (x\dashv y)\dashv z = x \dashv (y\dashv z), 
(x\dashv y)\dashv z = x \dashv (y\vdash z), 
(x\dashv y)\dashv z = x \dashv (x\circ z),\! {\rm\ (by\ (1))}\\
&&\hspace{-.7cm} (x\vdash y) \dashv z=x \vdash (y\dashv z), {\rm\ (by\ (2))}\\
&&\hspace{-.7cm} (x \vdash y)\vdash z=x\vdash (y\vdash z), (x\dashv y)\vdash z=x\vdash (y\vdash z),
(x\circ y)\vdash z= x \vdash (y\vdash z),\! {\rm\ (by\ (3))}\\
&&\hspace{-.7cm} \theta = \eta, {\rm\ with\ } \theta\in \{(x\dashv y)\circ z, (x\vdash y)\circ z, 
(x\circ y)\circ z, (x\circ y)\dashv z\}, \\
&&\hspace{-.7cm}\ \ \ \ \ \ \ \eta\in \{ x\vdash (y\circ z), x\circ (y\dashv z), x\circ (y\vdash z), x\circ (y\circ z)\},
{\rm\ (by\ (4)).}
\end{eqnarray*}
These are simplified (say by using Maple) to the following linearly independent subset.
\begin{eqnarray*}
&&(x\dashv y)\dashv z = x \dashv (y\dashv z), 
(x\dashv y)\dashv z = x \dashv (y\vdash z), 
(x\dashv y)\dashv z = x \dashv (x\circ z), \\
&& (x\vdash y) \dashv z=x \vdash (y\dashv z), \\
&&(x \vdash y)\vdash z=x\vdash (y\vdash z), (x\dashv y)\vdash z=x\vdash (y\vdash z),
(x\circ y)\vdash z= x \vdash (y\vdash z), \\
&& \theta = x\dashv (y\circ z) {\rm\ with\ } \theta\in \{(x\dashv y)\circ z, (x\vdash y)\circ z, 
(x\circ y)\circ z, (x\circ y)\dashv z\},\\
&& (x\circ y) \dashv z= \eta {\rm\ with\ } 
\eta\in \{ x\circ (y\dashv z), x\circ (y\vdash z), x\circ (y\circ z)\}.
\end{eqnarray*}
Since it has 14 elements, $\dfrelb_N$ has 4 elements and the dimension of 
$\dfgen_N^{\otimes 2}\oplus \dfgen_N^{\otimes 2}$ is 18, we see that this subset 
is a basis of $\dfrel_N^{\bot}$ and gives the relations of the dual operad of 
the NS operad. 
Of course, there are other choices for the basis. In particular, the relation 
$(x\circ y)\circ z = x\circ (y\circ z)$ is in $\dfrelb_N^{\bot}$. 
We will call this the {\bf associative Nijenhuis trialgebra}. It is similar 
to the 
associative trialgebra in that all of the three binary operations $\{\dashv,\vdash,\circ\}$
are associative. 

We remark that the products $\dftimes$ and $\dfotimes$ are not related by 
taking the dual, in contrast to the duality between the products $\bullet$ and
$\circ$ in~\cite{G-K}. For this we show that the dual of the quadri-algebra 
$(\dfgen_Q,\dfrel_Q)=(\dfgen_D,\dfrel_D) \dftimes (\dfgen_D,\dfrel_D)$ in 
Proposition~\ref{pp:quadri} is not 
$$(\dfgen_D,\dfrel_D)^! \dfotimes (\dfgen_D,\dfrel_D)^! 
    =(\dfgen_{AD},\dfrel_{AD}) \dfotimes (\dfgen_{AD},\dfrel_{AD}).$$
By definition, the dual of $(\dfgen_Q,\dfrel_Q)$ is given by 
$(\dfgen_{AQ},\dfrel_{AQ})$. Here 
$\dfgen_{AQ}$ has basis
$$ \Big\{ \dfpair{\dashv}{\dashv}, \dfpair{\dashv}{\vdash}, 
    \dfpair{\vdash}{\dashv}, \dfpair{\vdash}{\vdash} \Big\}$$ which is dual 
to the basis 
$$\left \{ \dfpair{\prec}{\prec},\dfpair{\prec}{\succ},
\dfpair{\succ}{\prec}, \dfpair{\succ}{\succ} \right \}$$ of $\dfgen_Q$. 
Also $\dfrel_{AQ}=\dfrel_{Q}^\perp$. 
By Eq. (\ref{eq:asdi}), $(\vdash\ot \dashv, \vdash \ot \dashv)$ is in 
$\dfrel_D^\perp=\dfrel_{AD}$. However, 
$\Big(\dfpair{\vdash}{\vdash}\ot \dfpair{\dashv}{\vdash},\dfpair{\vdash}{\vdash}
\ot \dfpair{\dashv}{\dashv}\Big )$ is not in $\dfrel_Q^\perp$. For example, 
$\Big (\dfpair{\succ}{\succ} \ot \dfpair{\prec}{\prec}, 
    \dfpair{\succ}{\succ}\ot \dfpair{\prec}{\prec} \Big )$ is in $\dfrel_Q$. 
But by Eq. (\ref{eq:pair1}) and (\ref{eq:pair2}), we have
\begin{eqnarray*}
\lefteqn{\Big \langle \Big (\dfpair{\succ}{\succ} \ot \dfpair{\prec}{\prec}, 
    \dfpair{\succ}{\succ}\ot \dfpair{\prec}{\prec} \Big ),
\Big(\dfpair{\vdash}{\vdash}\ot \dfpair{\dashv}{\vdash},
    \dfpair{\vdash}{\vdash}\ot \dfpair{\dashv}{\dashv}\Big ) 
           \Big \rangle_{2\dfgen_Q^{\ot 2}}  } \\
&=& \Big \langle \dfpair{\succ}{\succ} \ot \dfpair{\prec}{\prec}, 
\dfpair{\vdash}{\vdash}\ot \dfpair{\dashv}{\vdash} \Big \rangle_{\dfgen_Q^{\ot 2}}
- \Big \langle \dfpair{\succ}{\succ}\ot \dfpair{\prec}{\prec}, 
    \dfpair{\vdash}{\vdash}\ot \dfpair{\dashv}{\dashv}
           \Big \rangle_{\dfgen_Q^{\ot 2}}  \\
&=& \Big \langle \dfpair{\succ}{\succ}, \dfpair{\vdash}{\vdash} \Big\rangle_{\dfgen_Q}
\Big\langle \dfpair{\prec}{\prec}, \dfpair{\dashv}{\vdash} \Big\rangle_{\dfgen_Q}
- \Big \langle \dfpair{\succ}{\succ},\dfpair{\vdash}{\vdash} \Big \rangle_{\dfgen_Q}
\Big \langle \dfpair{\prec}{\prec}, \dfpair{\dashv}{\dashv} \Big \rangle_{\dfgen_Q} \\
&=& -1 
\end{eqnarray*}
So 
$\Big(\dfpair{\vdash}{\vdash}\ot \dfpair{\dashv}{\vdash},
    \dfpair{\vdash}{\vdash}\ot \dfpair{\dashv}{\dashv}\Big ) $ is not 
in $\dfrel_Q^\perp=\dfrel_{AQ}$.
Thus $\dfrel_{AD} \ot \dfgen_{AD}^{\ot 2} \not\subseteq \dfrel_{AQ}$ and 
therefore 
$\dfrel_{AD} \dfotimes \dfrel_{AD} \not\subseteq \dfrel_{AQ}$.

\section{\Loday algebras from linear operators}
\mlabel{sec:op}

\subsection{Operators on \Loday algebras}
A common method used to obtain a new operad structure from
a known operad structure is by means of a linear operator
on the known operad. Such linear operators include
(left, right and two-sided)
Rota-Baxter operators and Nijenhuis operators. Examples of such constructions
can be found in \cite{A-L,L-R2,Le1,Le2,Le3}. The constructions were usually
verified by checking the relations for each case of operators and operads.
We will verify the construction for (left, right and two-sided)
Rota-Baxter and Nijenhuis operators on all types of \Loday algebras.
We first give the definitions.

\begin{defn}
Let $D$ be an \Loday algebra of type $(\dfgen,\dfrel)$. A linear
operator $P$ on $D$ is called a {\bf(}two-sided{\bf)} {\bf
Rota-Baxter operator of weight $\lambda$} {\rm (resp. {\bf left
Rota-Baxter}, resp. {\bf right Rota-Baxter})} if, for each
$\dfop\in \dfgen$, we have
\begin{equation}
P(x)\dfop P(y)
=P(P(x)\dfop y+x\dfop P(y)+\lambda x\dfop y),\, x,\ y\in D$$
$${\rm (resp.\ }
P(x)\dfop P(y) =P(x\dfop P(y)),\, x,\ y\in D {\rm
)}\label{RBR}
\end{equation}
$${\rm (resp.\ }
P(x)\dfop P(y) =P(P(x)\dfop y),\, x,\ y\in D {)}.$$ A linear
operator $N$ on $D$ is called a {\bf Nijenhuis operator} if, for
each $\dfop\in \dfgen$, we have
\begin{equation}
N(x)\dfop N(y)
=N(N(x)\dfop y+x\dfop N(y) - N(x\dfop y)),\, x,\ y\in
D.\label{NR}
\end{equation}
\end{defn}

\begin{theorem}
Let $D$ be an \Loday algebra of type
$(\dfgen,\dfrel)$.
\begin{enumerate}
\item
Let $P$ be a Rota-Baxter operator of weight $\lambda$ on the \Loday algebra $D$.
For each $\dfop\in \dfgen$, define binary operations on $D$ by
$$ x\dfpair{\dfop}{\prec} y=x\dfop P(y), \
x\dfpair{\dfop}{\succ} y = P(x) \dfop y,\
x \dfpair{\dfop}{\circ} y= \lambda x \dfop y.$$
Then these operations define an \Loday algebra on $D$
of type $(\dfgen\otimes \dfgen_T, \dfrel\dftimes \dfrel_T)$
where $(\dfgen_T,\dfrel_T)$ is the type of dendriform
trialgebra.
\item
Let $N$ be a Nijenhuis operator on the \Loday algebra $D$. For each $\dfop\in \dfgen$,
define binary operations on $D$ by
$$ x\dfpair{\dfop}{\prec} y=x\dfop N(y), \
x\dfpair{\dfop}{\succ} y = N(x) \dfop y,\ x
\dfpair{\dfop}{\bullet} y= -N (x \dfop y).$$ Then these operations
define an \Loday algebra on $D$ of type $(\dfgen\otimes \dfgen_N,
\dfrel\dftimes \dfrel_N)$ where $(\dfgen_N,\dfrel_N)$ is the type
of NS algebra in Eq. (\ref{eq:NS}).
\item
Let $P$ be a left (resp. right) Rota-Baxter operator on the dendriform
algebra $D$. For each $\dfop\in \dfgen$, define binary operations
on $D$ by
$$ x\dfpair{\dfop}{\succ} y = P(x) \dfop y,
x \dfpair{\dfop}{\stackrel{\rightarrow}{\dft}} y=  x \dfop P(y),$$
$${\rm (resp.\ } x\dfpair{\dfop}{\prec} y=x\dfop P(y),
x \dfpair{\dfop}{\stackrel{\leftarrow}{\dft}} y= P(x) \dfop y {\rm
)}.$$ Then these operations define an \Loday algebra on $D$ of type
$(\dfgen\otimes \dfgen_L, \dfrel\dftimes \dfrel_L)$ {(}(resp.
$(\dfgen\otimes \dfgen_R, \dfrel\dftimes \dfrel_R)${\rm)} where
$(\dfgen_L,\dfrel_L)$ {\rm(}resp. $(\dfgen_R,\dfrel_R)${\rm)} is
the type of $L$-dipterous) {\rm (}resp. $L$-anti-dipterous{\rm)}
algebra in Eq. (\ref{eq:dip}) (resp. in Eq. (\ref{eq:anti})).
\end{enumerate}
\label{thm:op}
\end{theorem}

\begin{proof}
(1). Recall that
\allowdisplaybreaks{\begin{eqnarray*} \dfrelb_T=&
\{(\prec\otimes\prec,\prec\otimes\star),
(\succ\otimes\prec,\succ\otimes\prec), (\star\otimes
\succ,\succ\otimes\succ),
(\succ\otimes\circ,\succ\otimes\circ), \\
& (\prec\otimes\circ,\circ\otimes\succ),(\circ\otimes\prec,\circ\otimes\prec),
(\circ\otimes\circ,\circ\otimes\circ)\}.
\end{eqnarray*}}
To prove our claim we only need to prove that, for each
$$f=\sum_{j=1}^r(\dfoa_j\otimes\dfob_j,
\dfoc_j\otimes\dfod_j)\in \dfrel$$
and each of the 7 pairs in $\dfrelb_T$, the ``tensor product" of
the two is a relation on $\dfgen\otimes\dfgen_T$.
We consider the 7 cases separately.

Case 1. The pair is $(\prec\otimes\prec,\prec\otimes\star)$.
We need to verify that
$$f\dftimes (\prec\otimes \prec, \prec\otimes \star)
= \sum_{j=1}^{r}\left(\dfpair{\dfoa_j}{\prec}\otimes\dfpair{\dfob_j}{\prec},
\dfpair{\dfoc_j}{\prec}\otimes\dfpair{\dfod_j}{\star}\right)$$
is a relation on $D$. We have

\allowdisplaybreaks{\begin{eqnarray*} \lefteqn{\sum_{j=1}^{r}
\left (x \dfpair{\dfoa_j}{\prec}y \right )
\dfpair{\dfob_j}{\prec}z
= \sum_{j=1}^{r}( x \dfoa_j P(y)) \dfob_j P(z)} \\
&=& \sum_{j=1}^{r}  x \dfoc_j (P(y) \dfod_j P(z))\\
&=& \sum_{j=1}^{r}  x \dfoc_j P( P(y)\dfod_j z + y\dfod_j P(z) +
    \lambda y \dfod_j z)\\
&=& \sum_{j=1}^r  x \dfpair{\dfoc_j}{\prec} (P(y)\dfod_j z
    + y\dfod_j P(z) +  \lambda y \dfod_j z)\\
&=& \sum_{j=1}^r  x \dfpair{\dfoc_j}{\prec} \left (y\left (\dfpair{\dfod_j}{\prec}
+\dfpair{\dfod_j}{\succ}+\dfpair{\dfod_j}{\circ}\right )z\right )\\
&=& \sum_{j=1}^r  x \dfpair{\dfoc_j}{\prec}
\left (y\dfpair{\dfod_j}{\star} z\right),
\end{eqnarray*}}
as is desired.

Case 2. The pair is $(\succ\otimes\prec,\succ\otimes\prec)$. Then
we have
\allowdisplaybreaks{\begin{eqnarray*}
\lefteqn{\sum_{j=1}^{r} \left (x \dfpair{\dfoa_j}{\succ}y\right )
\dfpair{\dfob_j}{\prec}z
= \sum_{j=1}^{r} (P(x) \dfoa_j y) \dfob_j P(z)} \\
&=& \sum_{j=1}^{r}  P(x) \dfoc_j (y \dfod_j P(z))\\
&=& \sum_{j=1}^r  x \dfpair{\dfoc_j}{\succ} \left (y \dfpair{\dfod_j}{\prec} z\right).
\end{eqnarray*}}

Case 3. The pair is $(\star\otimes \succ,\succ\otimes\succ)$. We
have
\allowdisplaybreaks{\begin{eqnarray*} \sum_{j=1}^{r} \left (x
\dfpair{\dfoa_j}{\star}y \right )\dfpair{\dfob_j}{\succ}z &=&
\sum_{j=1}^{r} P(x\dfoa_j P(y) +P(x)\dfoa_j y+\lambda x\dfoa_j y)
\dfob_j z\\
&=& \sum_{j=1}^{r} (P(x)\dfoa_j P(y))\dfob_j z\\
&=& \sum_{j=1}^{r} P(x) \dfoc_j (P(y) \dfod_j z)\\
&=& \sum_{j=1}^{r} x \dfpair{\dfoc_j}{\succ} \left (y \dfpair{\dfod_j}{\succ} z\right).
\end{eqnarray*}}

Case 4. The pair is $(\succ\otimes\circ,\succ\otimes\circ)$. We
have
\allowdisplaybreaks{\begin{eqnarray*} \lefteqn{
\sum_{j=1}^{r}\left (x \dfpair{\dfoa_j}{\succ}y\right)\dfpair{\dfob_j}{\circ}z } \\
&=& \sum_{j=1}^{r} (P(x) \dfoa_j y) \dfob_j \lambda z \\
&=& \sum_{j=1}^{r} P(x) \dfoc_j (y \dfod_j \lambda z) \\
&=& \sum_{j=1}^{r} x \dfpair{\dfoc_j}{\succ}
\left (y \dfpair{\dfod_j}{\circ} z\right).
\end{eqnarray*}}

The cases for the pairs
$(\prec\otimes\circ,\circ\otimes\succ)$,
$(\circ\otimes\prec,\circ\otimes\prec)$
and
$(\circ\otimes\circ,\circ\otimes\circ)$
are proved in the same way as Case 4.

\medskip
\noindent
2. Now we consider a Nijenhuis operator $N$.
The relation set of the NS-algebra is
\begin{eqnarray*}
\dfrelb_N=& \{(\prec\otimes\prec,\prec\otimes\star),
(\succ\otimes\prec,\succ\otimes\prec),
(\star\otimes \succ,\succ\otimes\succ),\\
& (\star\otimes\bullet+\bullet\otimes\prec,
\succ\otimes\bullet+\bullet\otimes \star)\}.
\end{eqnarray*}
So we only need to prove that, for each
$$f=\sum_{j=1}^{r}(\dfoa_j\otimes\dfob_j,
\dfoc_j\otimes\dfod_j)\in \dfrel$$
and each of the 4 pairs in $\dfrelb_N$, the ``tensor product" of
the two is a relation for $D$.
The verification of the first three pairs is the same as the first
three cases in the Rota-Baxter operator case. For the fourth case, taking
$f$ as above, we have
\begin{eqnarray*}
&&\sum_{j=1}^r \left (\left (x \dfpair{\dfoa_j}{\star}y\right )\dfpair{\dfob_j}{\bullet}z
+\left (x \dfpair{\dfoa_j}{\bullet}y \right )\dfpair{\dfob_j}{\prec}z\right ) 
\\
&=&\!\!\sum_{j=1}^r\!\! \left (- N((N(x)\dfoa_j y\!+\!\dfoa_j N(y)-N(x\dfoa_j y))\dfob_j z)
-N(x \dfoa_j y) \dfob_j N(z) \right )
\\
&=&\sum_{j=1}^r \Big(- N((N(x)\dfoa_j y+x\dfoa_j N(y)-N(x\dfoa_j y))\dfob_j z)
\\
&&-N(N(x \dfoa_j y) \dfob_j z+(x\dfoa_j y) \dfob_j N(z)
    -N((x\dfoa_j y)\dfob_j z))\Big )\\
&=& \sum_{j=1}^r \Big (- N((N(x)\dfoa_j y)\dfob_j z
+(x\dfoa_j N(y))\dfob_j z\\
&&+(x\dfoa_j y) \dfob_j N(z)
    -N((x\dfoa_j y)\dfob_j z))\Big )\\
&=& -N \left (\sum_{j=1}^r  (N(x)\dfoa_j y)\dfob_j z
+\sum_{j=1}^r  (x\dfoa_j N(y))\dfob_j z \right .\\
&& +\sum_{j=1}^r \left . (x\dfoa_j y) \dfob_j N(z)
    -N\left (\sum_{j=1}^r (x\dfoa_j y)\dfob_j z\right )\right )\\
&=& -N \left (\sum_{j=1}^r  N(x)\dfoc_j (y\dfod_j z )
+\sum_{j=1}^r  x\dfoc_j (N(y)\dfod_j z )\right .\\
&& +\sum_{j=1}^r \left .  x\dfoc_j (y \dfod_j N(z))
    -N\left (\sum_{j=1}^r  x\dfoc_j (y\dfod_j z)\right )\right ).
\end{eqnarray*}
On the other hand,
\allowdisplaybreaks{
\begin{eqnarray*} 
&&\sum_{j=1}^r \left( x \dfpair{\dfoc_j}{\succ}\left (
y\dfpair{\dfod_j}{\bullet}z\right)
+x \dfpair{\dfoc_j}{\bullet}\left (y\dfpair{\dfod_j}{\star}z \right)\right) \\
&=&\!\!\sum_{j=1}^r\!\! \left (- N(x) \dfoc_j\! N(y\dfod_j\! z)
\!-\! N(x \dfoc_j\! (N(y)\dfod_j z\! +\!y\dfod_j\! N(z)\! -\!N(y\dfod_j\! z))) \right )\\
&=& \sum_{j=1}^r \Big( - N(N(x) \dfoc_j (y\dfod_j z)+x \dfoc_j N(y\dfod_j z)
-N(x\dfoc_j (y\dfod_j z))) \\
&&- N(x \dfoc_j (N(y)\dfod_j z +y\dfod_j N(z) -N(y\dfod_j z)))\Big)\\
&=& \sum_{j=1}^r - N (N(x) \dfoc_j (y\dfod_j z)
-N(x\dfoc_j (y\dfod_j z)) \\
&&+ x \dfoc_j (N(y)\dfod_j z) +x\dfoc_j (y\dfod_j N(z)) ).
\end{eqnarray*}}
This verifies the last relation.

The proof of (3) is the same as (actually simpler than)
the proof of (1).
\end{proof}

\subsection{Algebras with commuting operators}
Using Theorem~\ref{thm:op} inductively, we get
\begin{coro}
Let $(D,\{P_i\}_i)$ be an algebra of type $(\dfgen,\dfrel)$ with
$k$ commuting linear operators $P_i$, such as Rota-Baxter, Nijenhuis, 
left Rota-Baxter or right Rota-Baxter operators.
We obtain on $D$ an algebra structure of type
$(\dfgen,\dfrel)\dftimes \left ({\dftimes}_{i}
(\dfgen_i,\dfrel_i)\right )$ where
$(\dfgen_i,\dfrel_i)$ is the algebra type corresponding to
the operator $P_i,\ 1\leq i\leq k$.
\mlabel{co:op}
\end{coro}

\subsection{Examples}
Some of the previous known constructions can be obtained as special cases
of Theorem~\ref{thm:op} and Corollary~\ref{co:op}.
They can be found
\begin{enumerate}
\item
in \cite{Ag} where a Rota-Baxter operator of weight zero on an associative
algebra is used to construct a dendriform dialgebra;
\item
in \cite{EF1,Le1} where a Rota-Baxter operator of non-zero weight on an
associative algebra is used to construct a dendriform trialgebra;
\item
in \cite{A-L} where a Rota-Baxter operator of weight zero on a dendriform
dialgebra is used to construct a dendriform quadri-algebra;
\item
in \cite{A-L} where a pair of commuting Rota-Baxter operator of weight zero
on an associative algebra is used to construct a dendriform
quadri-algebra;
\item
in \cite{Le1} where a Rota-Baxter operator of non-zero weight on a
dendriform trialgebra is used to construct an ennea-algebra;
\item
in \cite{Le1} where a pair of commuting Rota-Baxter operators of non-zero
weight on an associative algebra is used to construct an ennea-algebra;
\item
in \cite{Le2} where a Nijenhuis operator on an associative algebra
is used to construct a NS-algebra;
\item
in \cite{Le2} where a Nijenhuis operator on a dendriform trialgebra
is used to construct a dendriform-Nijenhuis algebra.
\item
in \cite{Le3} where a left (resp. right) Rota-Baxter operator on an
associative algebra is used to construct a $L$-dipterous
(an anti-$L$-dipterous) algebra;
\item
in \cite{Le3} where a commuting pair of a right and a left Rota-Baxter
operators (resp. of two left Rota-Baxter operators) on an associative algebra
is used to construct an $M_1$ algebra (resp. an $M_2$ algebra);
\item
in \cite{Le3} where a set of three pairwise commuting Rota-Baxter operators
of weight zero on an associative algebra is used construct
an octo-algebra.
\end{enumerate}

\medskip
We end this paper with two more examples. It is easy to verify
that if $P$ is a Rota-Baxter operator of weight $\lambda \in \bfk$
on an \Loday algebra $D$ of type $(\dfgen,\dfrel)$, i.e. $P$
satisfies relation (\ref{RBR}), then the operator
$\tilde{P}:=-\lambda\, \id_{D}-P$ also is a Rota-Baxter operator
of weight $\lambda$ on $D$. Thus $D$ has two commuting Rota-Baxter
operators $P$ and $\tilde{P}$ of weight $\lambda$. So by
Corollary~\ref{co:op}, $D$ is equipped with an \Loday algebra of
type
$$(\dfgen,\dfrel)\dftimes (\dfgen_T,\dfrel_T) \dftimes (\dfgen_T,\dfrel_T)
=(\dfgen,\dfrel) \dftimes (\dfgen_E,\dfrel_E).$$
Here $(\dfgen_T,\dfrel_T)$ is the type of dendriform trialgebra and
$(\dfgen_E,\dfrel_E)$ is the type of dendriform ennea-algebra.

Similarly, if $N$ is a Nijenhuis operator on $D$, i.e. it
satisfies relation (\ref{NR}), then so is the operator
$\tilde{N}=id_D-N$. Thus $D$ is equipped with an \Loday algebra of
type $(\dfgen,\dfrel)\dftimes (\dfgen_N,\dfrel_N) \dftimes
(\dfgen_N,\dfrel_N).$ Here $(\dfgen_N,\dfrel_N)$ is the type of NS
algebra.

Returning to the case when $P$ is a Rota-Baxter operator,
define for each $\dfop\in \dfgen$, binary operations on $D$ by
$$ x\dfpair{\dfop}{\prec} y=x\dfop P(y), \
x\dfpair{\dfop}{\succ} y = -\tilde{P}(x) \dfop y.$$ Then these
operations define an \Loday algebra on $D$ of type $(\dfgen\otimes
\dfgen_D, \dfrel\dftimes \dfrel_D)$ where $(\dfgen_D,\dfrel_D)$ is
the type of dendriform dialgebra. The proof is the same as the one given
in~\cite{EF1}.


{\em Acknowledgements}: The first author would like to thank the
Ev. Studienwerk Villigst and the theory department of the
Physikalisches Institut, at Bonn University for generous support.
We are grateful to Prof. J.-L. Loday for suggestions on an earlier draft 
of this paper 
and for directing us to his recent paper~\cite{Lo7} where products 
and dual are also defined. We have modified some concepts and notations 
in our paper to be consistent with his.
We also thank the referee for helpful comments, 
especially for pointing out that the products defined in this
paper are different from the operad products defined in~\cite{G-K}, 
a fact that was also communicated to us by Loday. 


\vspace{-.5cm}

\addcontentsline{toc}{section}{\numberline {}References}

\end{document}